\documentclass[]{amsart}
\usepackage{amsmath, amsthm, amscd, amsfonts, amssymb, graphicx, color}
\usepackage[bookmarksnumbered, colorlinks, plainpages]{hyperref}

\theoremstyle{definition}

\theoremstyle{remark}

\numberwithin{equation}{section}

\begin{document}
\setcounter{page}{1}
\begin{center}
{\bf THE TOPOLOGICAL CENTERS    AND FACTORIZATION PROPERTIES OF MODULE ACTIONS AND $\ast-INVOLUTION$ ALGEBRAS}
\end{center}
\address{}
\address{}
\title[]{}
\dedicatory{}
\title[]{}
\author[]{KAZEM HAGHNEJAD AZAR  }

\address{}

\address{}

\dedicatory{}



\begin{abstract}
 For Banach left and right module actions, we  extend some propositions from Lau and $\ddot{U}lger$ into general situations  and we establish  the relationships between  topological centers of module actions.  We also introduce the new concepts as $Lw^*w$-property and $Rw^*w$-property for Banach $A-bimodule$ $B$ and we obtain some conclusions in the topological center of module actions and Arens regularity of Banach algebras. we also study some  factorization properties of left module actions and we find some relations of them and topological centers of module actions. For Banach algebra $A$, we extend the definition of  $\ast-involution$ algebra into Banach $A-bimodule$ $B$ with some results in the factorizations of $B^*$. We have some applications in  group algebras.
 \end{abstract} \maketitle

\begin{center}
{\bf  1.Introduction and Preliminaries}
\end{center}
\vspace{0.5cm}
\noindent As is well-known [1], the second dual $A^{**}$ of $A$ endowed with the either Arens multiplications is a Banach algebra. The constructions of the two Arens multiplications in $A^{**}$ lead us to definition of topological centers for $A^{**}$ with respect both Arens multiplications. The topological centers of Banach algebras, module actions and applications of them  were introduced and discussed in [6, 8, 13, 14, 15, 16, 17, 21, 22], and they have attracted by some attentions.

\noindent Now we introduce some notations and definitions that we used
throughout  this paper.\\
 Let $A$ be a Banach algebra. We
say that a  net $(e_{\alpha})_{{\alpha}\in I}$ in $A$ is a left
approximate identity $(=LAI)$ [resp. right
approximate identity $(=RAI)$] if,
 for each $a\in A$,   $e_{\alpha}a\longrightarrow a$ [resp. $ae_{\alpha}\longrightarrow a$]. For $a\in A$
 and $a^\prime\in A^*$, we denote by $a^\prime a$
 and $a a^\prime$ respectively, the functionals on $A^*$ defined by $\langle  a^\prime a,b\rangle  =\langle  a^\prime,ab\rangle  =a^\prime(ab)$ and $\langle  a a^\prime,b\rangle  =\langle  a^\prime,ba\rangle  =a^\prime(ba)$ for all $b\in A$.
  The Banach algebra $A$ is embedded in its second dual via the identification
 $\langle  a,a^\prime\rangle  $ - $\langle  a^\prime,a\rangle  $ for every $a\in
A$ and $a^\prime\in
A^*$. We denote the set   $\{a^\prime a:~a\in A~ and ~a^\prime\in
  A^*\}$ and
  $\{a a^\prime:~a\in A ~and ~a^\prime\in A^*\}$ by $A^*A$ and $AA^*$, respectively, clearly these two sets are subsets of $A^*$.  Let $A$ has a $BAI$. If the
equality $A^*A=A^*,~~(AA^*=A^*)$ holds, then we say that $A^*$
factors on the left (right). If both equalities $A^*A=AA^*=A^*$
hold, then we say
that $A^*$  factors on both sides.
 Let $X,Y,Z$ be normed spaces and $m:X\times Y\rightarrow Z$ be a bounded bilinear mapping. Arens in [1] offers two natural extensions $m^{***}$ and $m^{t***t}$ of $m$ from $X^{**}\times Y^{**}$ into $Z^{**}$ as following:\\
 \noindent1. $m^*:Z^*\times X\rightarrow Y^*$,~~~~~given by~~~$\langle  m^*(z^\prime,x),y\rangle  =\langle  z^\prime, m(x,y)\rangle  $ ~where $x\in X$, $y\in Y$, $z^\prime\in Z^*$,\\
 2. $m^{**}:Y^{**}\times Z^{*}\rightarrow X^*$,~~given by $\langle  m^{**}(y^{\prime\prime},z^\prime),x\rangle  =\langle  y^{\prime\prime},m^*(z^\prime,x)\rangle  $ ~where $x\in X$, $y^{\prime\prime}\in Y^{**}$, $z^\prime\in Z^*$,\\
3. $m^{***}:X^{**}\times Y^{**}\rightarrow Z^{**}$,~ given by~ ~ ~$\langle  m^{***}(x^{\prime\prime},y^{\prime\prime}),z^\prime\rangle  $ $=\langle  x^{\prime\prime},m^{**}(y^{\prime\prime},z^\prime)\rangle  $ \\~where ~$x^{\prime\prime}\in X^{**}$, $y^{\prime\prime}\in Y^{**}$, $z^\prime\in Z^*$.\\
The mapping $m^{***}$ is the unique extension of $m$ such that $x^{\prime\prime}\rightarrow m^{***}(x^{\prime\prime},y^{\prime\prime})$ from $X^{**}$ into $Z^{**}$ is $weak^*-to-weak^*$ continuous for every $y^{\prime\prime}\in Y^{**}$, but the mapping $y^{\prime\prime}\rightarrow m^{***}(x^{\prime\prime},y^{\prime\prime})$ is not in general $weak^*-to-weak^*$ continuous from $Y^{**}$ into $Z^{**}$ unless $x^{\prime\prime}\in X$. Hence the first topological center of $m$ may  be defined as following
$$Z_1(m)=\{x^{\prime\prime}\in X^{**}:~~y^{\prime\prime}\rightarrow m^{***}(x^{\prime\prime},y^{\prime\prime})~~is~~weak^*-to-weak^*-continuous\}.$$
Let now $m^t:Y\times X\rightarrow Z$ be the transpose of $m$ defined by $m^t(y,x)=m(x,y)$ for every $x\in X$ and $y\in Y$. Then $m^t$ is a continuous bilinear map from $Y\times X$ to $Z$, and so it may be extended as above to $m^{t***}:Y^{**}\times X^{**}\rightarrow Z^{**}$.
 The mapping $m^{t***t}:X^{**}\times Y^{**}\rightarrow Z^{**}$ in general is not equal to $m^{***}$, see [1], if $m^{***}=m^{t***t}$, then $m$ is called Arens regular. The mapping $y^{\prime\prime}\rightarrow m^{t***t}(x^{\prime\prime},y^{\prime\prime})$ is $weak^*-to-weak^*$ continuous for every $y^{\prime\prime}\in Y^{**}$, but the mapping $x^{\prime\prime}\rightarrow m^{t***t}(x^{\prime\prime},y^{\prime\prime})$ from $X^{**}$ into $Z^{**}$ is not in general  $weak^*-to-weak^*$ continuous for every $y^{\prime\prime}\in Y^{**}$. So we define the second topological center of $m$ as
$$Z_2(m)=\{y^{\prime\prime}\in Y^{**}:~~x^{\prime\prime}\rightarrow m^{t***t}(x^{\prime\prime},y^{\prime\prime})~~is~~weak^*-to-weak^*-continuous\}.$$
It is clear that $m$ is Arens regular if and only if $Z_1(m)=X^{**}$ or $Z_2(m)=Y^{**}$. Arens regularity of $m$ is equivalent to the following
$$\lim_i\lim_j\langle  z^\prime,m(x_i,y_j)\rangle  =\lim_j\lim_i\langle  z^\prime,m(x_i,y_j)\rangle  ,$$
whenever both limits exist for all bounded sequences $(x_i)_i\subseteq X$ , $(y_i)_i\subseteq Y$ and $z^\prime\in Z^*$, see [6, 18].\\
 The regularity of a normed algebra $A$ is defined to be the regularity of its algebra multiplication when considered as a bilinear mapping. Let $a^{\prime\prime}$ and $b^{\prime\prime}$ be elements of $A^{**}$, the second dual of $A$. By $Goldstin^,s$ Theorem [6, P.424-425], there are nets $(a_{\alpha})_{\alpha}$ and $(b_{\beta})_{\beta}$ in $A$ such that $a^{\prime\prime}=weak^*-\lim_{\alpha}a_{\alpha}$ ~and~  $b^{\prime\prime}=weak^*-\lim_{\beta}b_{\beta}$. So it is easy to see that for all $a^\prime\in A^*$,
$$\lim_{\alpha}\lim_{\beta}\langle  a^\prime,m(a_{\alpha},b_{\beta})\rangle  =\langle  a^{\prime\prime}b^{\prime\prime},a^\prime\rangle  $$ and
$$\lim_{\beta}\lim_{\alpha}\langle  a^\prime,m(a_{\alpha},b_{\beta})\rangle  =\langle  a^{\prime\prime}ob^{\prime\prime},a^\prime\rangle  ,$$
where $a^{\prime\prime}b^{\prime\prime}$ and $a^{\prime\prime}ob^{\prime\prime}$ are the first and second Arens products of $A^{**}$, respectively, see [6, 14, 18].\\
The mapping $m$ is left strongly Arens irregular if $Z_1(m)=X$ and $m$ is right strongly Arens irregular if $Z_2(m)=Y$.\\
This paper is organized as follows.\\
{\bf a)} In section two, for  a  Banach $A-bimodule$, we have
 \begin{enumerate}
 \item $a^{\prime\prime}\in Z_{B^{**}}(A^{**})$ if and only if $\pi_\ell^{****}(b^\prime,a^{\prime\prime})\in B^*$ for all $b^\prime\in B^*$.

\item $F\in Z_{B^{**}}((A^*A)^*)$ if and only if $\pi_\ell^{****}(g,F)\in B^{*}$ for all $g\in B^{*}$.
\item  $G\in Z_{(A^*A)^*}(B^{**})$ if and only if $\pi_r^{****}(g,G)\in A^*A$ for all $g\in B^{*}$.
\item  Let $B$ has a $BAI$ $(e_\alpha)_\alpha\subseteq A$ such that
$e_\alpha \stackrel{w^*} {\rightarrow}e^{\prime\prime}$. Then if ${Z}^t_{e^{**}}(B^{**})=B^{**}$ [ resp. ${Z}_{e^{**}}(B^{**})=B^{**}$] and $B^*$ factors on the left [resp. right], but not on the right [resp. left], then ${Z}_{B^{**}}(A^{**})\neq {Z}^t_{B^{**}}(A^{**})$.

\item $B^*A\subseteq wap_\ell(B)$
if and only if  $AA^{**}\subseteq Z_{B^{**}}(A^{**})$.
\item  Let $b^\prime\in B^*$. Then $b^\prime\in wap_\ell(B)$ if and only if the adjoint of  the mapping $\pi_\ell^*(b^\prime,~):A\rightarrow B^*$ is $weak^*-to-weak$ continuous.

Then $(\alpha)\Rightarrow (\beta)\Rightarrow (\gamma)$.\\
If we take $T\in LM(A,B)$ and if $B$ has a sequential $BAI$ and it is $WSC$, then ($\alpha$), ($\beta$) and ($\gamma$) are equivalent.
\end{enumerate}
{\bf b)} In section three, for  a Banach  $A-bimodule$ $B$, we define $Left-weak^*-to-weak$ property [=$Rw^*w-$ property]  and $Right-weak^*-to-weak$ property [=$Rw^*w-$ property] for Banach algebra $A$ and we show that
\begin{enumerate}
\item If $A^{**}=a_0A^{**}$ [resp.  $A^{**}=A^{**}a_0$] for some $a_0\in A$ and $a_0$ has $Rw^*w-$ property [resp. $Lw^*w-$ property], then $Z_{B^{**}}(A^{**})=A^{**}$.
\item If $B^{**}=a_0B^{**}$ [resp.  $B^{**}=B^{**}a_0$] for some $a_0\in A$ and $a_0$ has $Rw^*w-$ property [resp. $Lw^*w-$ property] with respect to $B$, then $Z_{A^{**}}(B^{**})=B^{**}$.
\item  If $B^*$ factors on the left [resp. right] with respect to $A$ and $A$ has $Rw^*w-$ property [resp. $Lw^*w-$ property], then $Z_{B^{**}}(A^{**})=A^{**}$.
\item  If $B^*$ factors on the left [resp. right] with respect to $A$ and $A$ has $Rw^*w-$ property [resp. $Lw^*w-$ property] with respect $B$, then $Z_{A^{**}}(B^{**})=B^{**}$.
\item  If $a_0\in A$ has $Rw^*w-$ property with respect to $B$, then $a_0A^{**}\subseteq Z_{B^{**}}(A^{**})$ and $a_0B^*\subseteq wap_\ell(B)$.
 \item    Assume that $AB^*\subseteq wap_\ell B$. If
$B^*$ strong factors on the left [resp. right], then $A$ has $Lw^*w-$ property [resp. $Rw^*w-$ property ] with respect to $B$.
\item Assume that $AB^*\subseteq wap_\ell B$. If
$B^*$ strong factors on the left [resp. right], then $A$ has $Lw^*w-$ property [resp. $Rw^*w-$ property ] with respect to $B$.
\end{enumerate}
{\bf c)} In section four, for  a left Banach  $A-module$ $B$, we study some relationships between factorization properties of left module action and topological centers of it.
\begin{enumerate}
\item Let  $(e_{\alpha})_{\alpha}\subseteq A$ be a $LBAI$ for $B$. Then the following assertions hold.\\
 i) For each $b^\prime\in B^*$, $\pi_\ell^*(b^\prime,e_\alpha)\stackrel{w^*} {\rightarrow}b^\prime$.\\
 ii)  $B^*$ factors on the left with respect to $A$ if and only if $B^{**}$ has a $W^*LBAI$ $(e_{\alpha})_{\alpha}\subseteq A$.\\
iii)    $B^{**}$ has a $W^*LBAI$ $(e_{\alpha})_{\alpha}\subseteq A$ if and only if $B^{**}$ has a left unit element $e^{\prime \prime}\in A^{**}$ such that $e_{\alpha}\stackrel{w^*} {\rightarrow}e^{\prime \prime}$.\\
\item Suppose that $b^\prime\in wap_\ell(B)$. Let  $a^{\prime\prime}\in A^{**}$ and
$(a_{\alpha})_{\alpha}\subseteq A$ such that  $a_{\alpha} \stackrel{w^*} {\rightarrow}a^{\prime\prime}$ in $A^{**}$. Then we have
$$\pi^*_\ell(b^\prime,a_\alpha) \stackrel{w} {\rightarrow}\pi_\ell^{****}(b^\prime,a^{\prime\prime}).$$
\item  Let  $B^*$ factors on the left with respect to $A$. If $AA^{**}\subseteq  Z_{B^{**}}(A^{**})$, then  $Z_{B^{**}}(A^{**})=A^{**}$.
\item  Let $B^{**}$ has a $LBAI$ with respect to $A^{**}$. Then $B^{**}$ has a left unit with respect to $A^{**}$.
\item  Let $B$  has a $LBAI$ with respect to $A$. Then we have the following assertions.\\
i)   $B^*$ factors on the left with respect to $A$ if and only if for each $b^\prime\in B^*$, we have $\pi^*_\ell(b^\prime,e_\alpha) \stackrel{w} {\rightarrow}b^\prime$ in $B^*$.\\
ii)  $B$ factors on the left with respect to $A$ if and only if for each $b\in B$, we have $\pi^*_\ell(b,e_\alpha) \stackrel{w} {\rightarrow}b$ in $B$.
\item  Let $A$ has a $LBAI$ $(e_{\alpha})_{\alpha}\subseteq A$ such that $e_{\alpha} \stackrel{w^*} {\rightarrow}e^{\prime\prime}$ in $A^{**}$ where $e^{\prime\prime}$ is a left unit for $A^{**}$. Suppose that $Z^t_{e^{\prime\prime}}(B^{**})=B^{**}$. Then, $B$ factors on the right with respect to $A$ if and only if
$e^{\prime\prime}$ is a left unit for $B^{**}$.

\end{enumerate}
{\bf d)} In section five, for  a Banach  $A-bimodule$ $B$ we study $\ast-involution$ algebra on $B^{**}$ with respect to first Arens product with some results in the factorization in $B^*$, that is,   suppose that $(e_{\alpha})_{\alpha}\subseteq A$ is a $BAI$ for $B$. Then if $B^{**}$ is a Banach $\ast-involution$ algebra as $A^{**}$-module,  then $B^{**}$ is unital  as $A^{**}$-module and $B^*$ factors on the both side.\\\\\\

\begin{center}
\textbf{{ 2. The topological centers of module actions}}\\
\end{center}

\vspace{0.5cm}

 Let $B$ be a Banach $A-bimodule$, and let\\
$$\pi_\ell:~A\times B\rightarrow B~~~and~~~\pi_r:~B\times A\rightarrow B.$$
be the left and right module actions of $A$ on $B$. Then $B^{**}$ is a Banach $A^{**}-bimodule$ with module actions
$$\pi_\ell^{***}:~A^{**}\times B^{**}\rightarrow B^{**}~~~and~~~\pi_r^{***}:~B^{**}\times A^{**}\rightarrow B^{**}.$$
Similarly, $B^{**}$ is a Banach $A^{**}-bimodule$ with module actions\\
$$\pi_\ell^{t***t}:~A^{**}\times B^{**}\rightarrow B^{**}~~~and~~~\pi_r^{t***t}:~B^{**}\times A^{**}\rightarrow B^{**}.$$
We may therefore define the topological centers of the right and left module actions of $A$ on $B$ as follows:
$$Z_{A^{**}}(B^{**})=Z(\pi_r)=\{b^{\prime\prime}\in B^{**}:~the~map~~a^{\prime\prime}\rightarrow \pi_r^{***}(b^{\prime\prime}, a^{\prime\prime})~:~A^{**}\rightarrow B^{**}$$$$~is~~~weak^*-to-weak^*~continuous\}$$
$$Z_{B^{**}}(A^{**})=Z(\pi_\ell)=\{a^{\prime\prime}\in A^{**}:~the~map~~b^{\prime\prime}\rightarrow \pi_\ell^{***}(a^{\prime\prime}, b^{\prime\prime})~:~B^{**}\rightarrow B^{**}$$$$~is~~~weak^*-to-weak^*~continuous\}$$
$$Z_{A^{**}}^t(B^{**})=Z(\pi_\ell^t)=\{b^{\prime\prime}\in B^{**}:~the~map~~a^{\prime\prime}\rightarrow \pi_\ell^{t***}(b^{\prime\prime}, a^{\prime\prime})~:~A^{**}\rightarrow B^{**}$$$$~is~~~weak^*-to-weak^*~continuous\}$$
$$Z_{B^{**}}^t(A^{**})=Z(\pi_r^t)=\{a^{\prime\prime}\in A^{**}:~the~map~~b^{\prime\prime}\rightarrow \pi_r^{t***}(a^{\prime\prime}, b^{\prime\prime})~:~B^{**}\rightarrow B^{**}$$$$~is~~~weak^*-to-weak^*~continuous\}$$
We note also that if $B$ is a left(resp. right) Banach $A-module$ and $\pi_\ell:~A\times B\rightarrow B$~(resp. $\pi_r:~B\times A\rightarrow B$) is left (resp. right) module action of $A$ on $B$, then $B^*$ is a right (resp. left) Banach $A-module$. \\
We write $ab=\pi_\ell(a,b)$, $ba=\pi_r(b,a)$, $\pi_\ell(a_1a_2,b)=\pi_\ell(a_1,a_2b)$,  $\pi_r(b,a_1a_2)=\pi_r(ba_1,a_2)$,~\\
$\pi_\ell^*(a_1b^\prime, a_2)=\pi_\ell^*(b^\prime, a_2a_1)$,~
$\pi_r^*(b^\prime a, b)=\pi_r^*(b^\prime, ab)$,~ for all $a_1,a_2, a\in A$, $b\in B$ and  $b^\prime\in B^*$
when there is no confusion.\\\\

\noindent{\it{\bf Theorem 2-1.}} We have the following assertions.
\begin{enumerate}
\item  Assume that   $B$ is a Banach left $A-module$. Then, $a^{\prime\prime}\in Z_{B^{**}}(A^{**})$ if and only if $\pi_\ell^{****}(b^\prime,a^{\prime\prime})\in B^*$ for all $b^\prime\in B^*$.
\item  Assume that   $B$ is a Banach right $A-module$. Then, $b^{\prime\prime}\in Z_{A^{**}}(B^{**})$ if and only if $\pi_r^{****}(b^\prime,b^{\prime\prime})\in A^*$ for all $b^\prime\in B^*$.
\end{enumerate}
\begin{proof}
\begin{enumerate}
\item Let $b^{\prime\prime}\in B^{**}$. Then, for every $a^{\prime\prime}\in Z_{B^{**}}(A^{**})$, we have
$$\langle  \pi_\ell^{****}(b^\prime,a^{\prime\prime}),b^{\prime\prime}\rangle  =
\langle  b^\prime,\pi_\ell^{***}(a^{\prime\prime},b^{\prime\prime})\rangle  =
\langle  \pi_\ell^{***}(a^{\prime\prime},b^{\prime\prime}),b^\prime\rangle  $$
$$=\langle  \pi_\ell^{t***t}(a^{\prime\prime},b^{\prime\prime}),b^\prime\rangle  =
\langle  \pi_\ell^{t***}(b^{\prime\prime},a^{\prime\prime}),b^\prime\rangle  =
\langle  b^{\prime\prime},\pi_\ell^{t**}(a^{\prime\prime},b^\prime)\rangle  .$$
It follow that $\pi_\ell^{****}(b^\prime,a^{\prime\prime})=\pi_\ell^{t**}(a^{\prime\prime},b^\prime)\in B^*$.\\
Conversely,  let $a^{\prime\prime}\in A^{**}$ and let $\pi_\ell^{****}(a^{\prime\prime},b^\prime)\in B^*$ for all $b^{\prime}\in B^{*}$.  Let $(b_{\alpha}^{\prime\prime})_{\alpha}\subseteq B^{**}$ such that  $b^{\prime\prime}_{\alpha} \stackrel{w^*} {\rightarrow}b^{\prime\prime}$. Then for every $b^\prime\in B^*$, we have 
$$\langle  \pi_\ell^{***}(a^{\prime\prime},b_\alpha^{\prime\prime}),b^\prime\rangle  =
\langle  b^\prime,\pi_\ell^{***}(a^{\prime\prime},b_\alpha^{\prime\prime})\rangle  =\langle  \pi_\ell^{****}(b^\prime,a^{\prime\prime}),
b_\alpha^{\prime\prime}\rangle =\langle
b_\alpha^{\prime\prime},\pi_\ell^{****}(b^\prime,a^{\prime\prime})\rangle  $$
$$\rightarrow \langle  
b^{\prime\prime},\pi_\ell^{****}(b^\prime,a^{\prime\prime})\rangle  =\langle  \pi_\ell^{****}(b^\prime,a^{\prime\prime}),
b^{\prime\prime}\rangle=\langle  \pi_\ell^{***}(a^{\prime\prime},b^{\prime\prime}),b^\prime\rangle.$$ 
Consequently $a^{\prime\prime}\in Z_{B^{**}}(A^{**})$.
\item Proof is similar to (1).
\end{enumerate}\end{proof}

\vspace{0.4cm}
\noindent{\it{\bf Theorem 2-2.}} Assume that   $B$ is a Banach  $A-bimodule$. Then we have the following assertions.
\begin{enumerate}

\item $F\in Z_{B^{**}}((A^*A)^*)$ if and only if $\pi_\ell^{****}(g,F)\in B^{*}$ for all $g\in B^{*}$.
\item  $G\in Z_{(A^*A)^*}(B^{**})$ if and only if $\pi_r^{****}(g,G)\in A^*A$ for all $g\in B^{*}$.
\end{enumerate}
\begin{proof}
 \begin{enumerate}
 \item Let $F\in Z_{B^{**}}((A^*A)^*)$ and $(b_{\alpha}^{\prime\prime})_{\alpha}\subseteq B^{**}$ such that  $b^{\prime\prime}_{\alpha} \stackrel{w^*} {\rightarrow}b^{\prime\prime}$. Then for all $g\in B^{*}$, we have
$$\langle  \pi_\ell^{****}(g,F),b^{\prime\prime}_{\alpha}\rangle  =\langle  g,\pi_\ell^{***}(F,b^{\prime\prime}_{\alpha})\rangle  =
\langle  \pi_\ell^{***}(F,b^{\prime\prime}_{\alpha}),g\rangle  $$
$$\rightarrow \langle  \pi_\ell^{***}(F,b^{\prime\prime}),g\rangle  =\langle  \pi_\ell^{****}(g,F),b^{\prime\prime}\rangle  .$$
Thus, we conclude that $\pi_\ell^{****}(g,F)\in (B^{**},weak^*)^*=B^{*}$.\\
Conversely, let $\pi_\ell^{****}(g,F)\in B^{*}$ for $F\in (A^*A)^*$ and $g\in B^{*}$. Assume that $b^{\prime\prime}\in B^{**}$ and  $(b_{\alpha}^{\prime\prime})_{\alpha}\subseteq B^{**}$ such that  $b^{\prime\prime}_{\alpha} \stackrel{w^*} {\rightarrow}b^{\prime\prime}$. Then
$$\langle  \pi_\ell^{***}(F,b^{\prime\prime}_{\alpha}),g\rangle  =\langle  g,\pi_\ell^{***}(F,b^{\prime\prime}_{\alpha})\rangle
=\langle  \pi_\ell^{****}(g,F),b^{\prime\prime}_{\alpha}\rangle  $$
$$=\langle  b^{\prime\prime}_{\alpha},\pi_\ell^{****}(g,F)\rangle  \rightarrow \langle  b^{\prime\prime},\pi_\ell^{****}(g,F)\rangle  =\langle  \pi_\ell^{****}(g,F),b^{\prime\prime}\rangle  $$$$
=\langle  \pi_\ell^{***}(F,b^{\prime\prime}),g\rangle  .$$
It follow that $F\in Z_{B^{**}}((A^*A)^*)$.
\item Proof is similar to (1).
 \end{enumerate}
 \end{proof}

 \vspace{0.5cm}

In the preceding theorems, if we take $B=A$, we obtain some parts of Lemma 3.1 from [14].\\\\

\noindent An element $e^{\prime\prime}$ of $A^{**}$ is said to be a mixed unit if $e^{\prime\prime}$ is a
right unit for the first Arens multiplication and a left unit for
the second Arens multiplication. That is, $e^{\prime\prime}$ is a mixed unit if
and only if
for each $a^{\prime\prime}\in A^{**}$, $a^{\prime\prime}e^{\prime\prime}=e^{\prime\prime}o a^{\prime\prime}=a^{\prime\prime}$. By
[4, p.146], an element $e^{\prime\prime}$ of $A^{**}$  is  mixed
      unit if and only if it is a $weak^*$ cluster point of some BAI $(e_\alpha)_{\alpha \in I}$  in
      $A$.\\
Let $B$ be a Banach  $A-bimodule$ and $a^{\prime\prime}\in A^{**}$. We define the locally topological center of the left and right module actions of $a^{\prime\prime}$ on $B$, respectively, as follows\\
$$Z_{a^{\prime\prime}}^t(B^{**})=Z_{a^{\prime\prime}}^t(\pi_\ell^t)=\{b^{\prime\prime}\in B^{**}:~~~\pi^{t***t}_\ell(a^{\prime\prime},b^{\prime\prime})=
\pi^{***}_\ell(a^{\prime\prime},b^{\prime\prime})\},$$
$$Z_{a^{\prime\prime}}(B^{**})=Z_{a^{\prime\prime}}(\pi_r^t)=\{b^{\prime\prime}\in B^{**}:~~~\pi^{t***t}_r(b^{\prime\prime},a^{\prime\prime})=
\pi^{***}_r(b^{\prime\prime},a^{\prime\prime})\}.$$\\
Thus we have ~~~~~~~~~$$\bigcap_{a^{\prime\prime}\in A^{**}}Z_{a^{\prime\prime}}^t(B^{**})=Z_{A^{**}}^t(B^{**})=
Z(\pi_r^t),  $$ ~~~~~~~ ~~ $$~~~~~~~~~~~~~~~~~~~~~~~~~~~~\bigcap_{a^{\prime\prime}\in A^{**}}Z_{a^{\prime\prime}}(B^{**})=Z_{A^{**}}(B^{**})=
Z(\pi_r).$$

Let $B$ be a Banach $A-bimodule$. We say that $B$ is a left [resp. right] factors with respect to $A$, if $BA=B$ [resp. $AB=B$].\\\\

\noindent{\it{\bf Definition 2-3.}} Let $B$ be a Banach left  $A-module$ and $e^{\prime\prime}\in A^{**}$ be a mixed unit for $A^{**}$. We say that $e^{\prime\prime}$ is a left mixed unit for $B^{**}$, if
$$\pi_\ell^{***}(e^{\prime\prime},b^{\prime\prime})=\pi_\ell^{t***t}(e^{\prime\prime},b^{\prime\prime})=b^{\prime\prime},$$
for all $b^{\prime\prime}\in B^{**}$.\\
The definition of right mixed unit for $B^{**}$ is similar. $B^{**}$ has a mixed unit, if it has equal left and right mixed unit.\\
It is clear that if   $e^{\prime\prime}\in A^{**}$ is a left (resp. right) unit for $B^{**}$ and $Z_{e^{\prime\prime}}(B^{**})=B^{**}$, then $e^{\prime\prime}$ is left (resp. right) mixed unit for $B^{**}$.\\\\

\noindent{\it{\bf Theorem 2-4.}} Let $B$ be a   Banach $A-bimodule$ with a $BAI$ $(e_\alpha)_\alpha$ such that
$e_\alpha \stackrel{w^*} {\rightarrow}e^{\prime\prime}$. Then if ${Z}^t_{e^{\prime\prime}}(B^{**})=B^{**}$ [resp. ${Z}_{e^{\prime\prime}}(B^{**})=B^{**}$] and $B^*$ factors on the left [resp. right], but not on the right [resp. left], then ${Z}_{B^{**}}(A^{**})\neq {Z}^t_{B^{**}}(A^{**})$.
\begin{proof} Suppose that  $B^*$ factors on the left with respect to $A$, but not on the right. Let $(e_{\alpha})_{\alpha}\subseteq A$ be a BAI for $A$ such that  $e_{\alpha} \stackrel{w^*} {\rightarrow}e^{\prime\prime}$. Thus for all $b^{\prime}\in B^{*}$ there are $a\in A$ and $x^\prime\in B^*$ such that $x^\prime a=b^\prime$. Then for all $b^{\prime\prime}\in B^{**}$ we have
$$\langle  \pi_\ell^{***}(e^{\prime\prime},b^{\prime\prime}),b^\prime\rangle
=\langle  e^{\prime\prime},\pi_\ell^{**}(b^{\prime\prime},b^\prime)\rangle  =
\lim_\alpha\langle  \pi_\ell^{**}(b^{\prime\prime},b^\prime),e_{\alpha}\rangle  $$
$$=\lim_\alpha\langle  b^{\prime\prime},\pi_\ell^{*}(b^\prime,e_{\alpha})\rangle
=\lim_\alpha\langle  b^{\prime\prime},\pi_\ell^{*}(x^\prime a,e_{\alpha})\rangle  $$
$$=\lim_\alpha\langle  b^{\prime\prime},\pi_\ell^{*}(x^\prime ,ae_{\alpha})\rangle
=\lim_\alpha\langle  \pi_\ell^{**}(b^{\prime\prime},x^\prime) ,ae_{\alpha}\rangle  $$$$=\langle  \pi_\ell^{**}(b^{\prime\prime},x^\prime) ,a\rangle
=\langle  b^{\prime\prime},b^{\prime}\rangle  .$$
Thus $\pi^{***}_\ell(e^{\prime\prime},b^{\prime\prime})=b^{\prime\prime}$ consequently $B^{**}$ has left unit $A^{**}-module$. It follows that
$e^{\prime\prime}\in {Z}_{B^{**}}(A^{**})$.  If we take
${Z}_{B^{**}}(A^{**})= {Z}^t_{B^{**}}(A^{**})$, then $e^{\prime\prime}\in {Z}^t_{B^{**}}(A^{**})$. Then the mapping $b^{\prime\prime}\rightarrow\pi_r^{t***t}(b^{\prime\prime},e^{\prime\prime})$ is $weak^*-to-weak^*$ continuous from $B^{**}$ into $B^{**}$. Since $e_\alpha \stackrel{w^*} {\rightarrow}e^{\prime\prime}$,
$\pi_r^{t***t}(b^{\prime\prime},e_\alpha)\stackrel{w^*} {\rightarrow}\pi_r^{t***t}(b^{\prime\prime},e^{\prime\prime})$. Let $b^\prime\in B^*$ and $(b_\beta)_\beta\subseteq B$ such that $b_\beta \stackrel{w^*} {\rightarrow}b^{\prime\prime}$. Since ${Z}^t_{e^{\prime\prime}}(B^{**})=B^{**}$, we have the following quality
$$\langle  \pi_r^{t***t}(b^{\prime\prime},e^{\prime\prime}), b^\prime\rangle  =\lim_\alpha\langle  \pi_r^{t***t}(b^{\prime\prime},e_\alpha), b^\prime\rangle  =\lim_\alpha\langle  \pi_r^{t***}(e_\alpha,b^{\prime\prime}), b^\prime\rangle  $$$$=\lim_\alpha\lim_\beta\langle  \pi_r^{t***}(e_\alpha,b_\beta), b^\prime\rangle
=
\lim_\alpha\lim_\beta\langle  \pi_r^{}(b_\beta ,e_\alpha), b^\prime\rangle  $$$$=\lim_\alpha\lim_\beta\langle   b^\prime, \pi_r^{}(b_\beta ,e_\alpha)\rangle
=\lim_\beta\lim_\alpha\langle   b^\prime, \pi_r^{}(b_\beta ,e_\alpha)\rangle  $$$$=\lim_\beta\langle   b^\prime, b_\beta \rangle  =
\langle  b^{\prime\prime},b^\prime\rangle  .$$
 Thus $\pi_r^{t***t}(b^{\prime\prime},e^{\prime\prime})=\pi_r^{***}(b^{\prime\prime},e^{\prime\prime})=b^{\prime\prime}$.
 It follows that $B^{**}$ has a right unit.
 Suppose that $b^{\prime\prime}\in B^{**}$ and $(b_\beta)_\beta\subseteq B$ such that $b_\beta \stackrel{w^*} {\rightarrow}b^{\prime\prime}$.  Then for all $b^\prime\in B^*$ we have
$$\langle  b^{\prime\prime},b^{\prime}\rangle  =\langle  \pi_r^{***}(b^{\prime\prime},e^{\prime\prime}),b^\prime\rangle  =
\langle  b^{\prime\prime},\pi_r^{**}(e^{\prime\prime},b^\prime)\rangle  =\lim_\beta\langle  \pi_r^{**}(e^{\prime\prime},b^\prime),b_\beta\rangle  $$
$$=\lim_\beta \langle  e^{\prime\prime},\pi_r^{*}(b^\prime,b_\beta)\rangle
=\lim_\beta \lim_\alpha\langle  \pi_r^{*}(b^\prime,b_\beta),e_\alpha\rangle  $$
$$=\lim_\beta \lim_\alpha\langle  \pi_r^{*}(b^\prime,b_\beta),e_\alpha\rangle  =\lim_\beta \lim_\alpha\langle  b^\prime,\pi_r(b_\beta,e_\alpha)\rangle  $$
$$=\lim_\alpha \lim_\beta\langle  \pi_r^{***}(b_\beta,e_\alpha),b^\prime\rangle
=\lim_\alpha \lim_\beta\langle  b_\beta,\pi_r^{**}(e_\alpha,b^\prime)\rangle  $$
$$= \lim_\alpha\langle  b^{\prime\prime},\pi_r^{**}(e_\alpha,b^\prime)\rangle  .$$
It follows that $weak-\lim_\alpha\pi_r^{**}(e_\alpha,b^\prime)=b^\prime$. So by Cohen factorization theorem, $B^*$ factors on the right that is contradiction.\end{proof}
\vspace{0.5cm}

\noindent{\it{\bf Corollary 2-5.}} Let $B$ be a   Banach $A-bimodule$ and $e^{\prime\prime}\in A^{**}$ be a  mixed unit for $B^{**}$. If $B^*$ factors on the left, but not on the right, then ${Z}_{B^{**}}(A^{**})\neq {Z}^t_{B^{**}}(A^{**})$.\\\\

For a Banach algebra $A$, if $e^{\prime\prime}$ is a left mixed unit for $A^{**}$, then it is clear that ${Z}^t_{e^{\prime\prime}}(A^{**})=A^{**}$. Then by using preceding theorem, we obtain Proposition 2.10 from [14].\\
\vspace{0.4cm}
\noindent We say that a Banach space $B$ is weakly complete, if for every  $(b_{\alpha})_{\alpha} \subseteq B$, we have  $b_{\alpha} \stackrel{w} {\rightarrow}b^{\prime\prime}$ in  $B^{**}$, it follows that $b^{\prime\prime}\in B$.\\\\

\noindent{\it{\bf Theorem 2-6.}} Suppose that $B$ is a weakly complete Banach space. Then we have the following assertions.
\begin{enumerate}
\item Let $B$ be a Banach left $A-module$ and  $B^{**}$ has a left mixed unit $e^{\prime\prime}\in A^{**}$. If $AB^{**}\subseteq B$, then $B$ is reflexive.
\item Let $B$ be a Banach right $A-module$ and $B^{**}$ has a right mixed unit $e^{\prime\prime}\in A^{**}$.  If $Z_{A^{**}}(B^{**})A\subseteq B$, then $Z_{A^{**}}(B^{**})= B$.
\end{enumerate}
\begin{proof}
\begin{enumerate}
\item Assume that $b^{\prime\prime}\in B^{**}$. Since $e^{\prime\prime}$ is also mixed unit for $A^{**}$, there is a $BAI$  $(e_{\alpha})_{\alpha}\subseteq A$ for $A$ such that $e_{\alpha} \stackrel{w^*} {\rightarrow}e^{\prime\prime}$. Then
$\pi_\ell^{***}(e_{\alpha},b^{\prime\prime}) \stackrel{w^*} {\rightarrow}\pi_\ell^{***}(e^{\prime\prime},b^{\prime\prime})=b^{\prime\prime}$ in $B^{**}$. Since $AB^{**}\subseteq B$, we have
$\pi_\ell^{***}(e_{\alpha},b^{\prime\prime})\in B$.  Consequently $\pi_\ell^{***}(e_{\alpha},b^{\prime\prime}) \stackrel{w} {\rightarrow}\pi_\ell^{***}(e^{\prime\prime},b^{\prime\prime})=b^{\prime\prime}$ in $B$.
Since $B$ is a weakly complete, $b^{\prime\prime}\in B$, and so $B$ is reflexive.
\item  Since  $b^{\prime\prime}\in Z_{A^{**}}(B^{**})$, we have  $\pi_r^{***}(b^{\prime\prime},e_{\alpha}) \stackrel{w^*} {\rightarrow}\pi_r^{***}(b^{\prime\prime},e^{\prime\prime})=b^{\prime\prime}$ in $B^{**}$.
    Since $Z_{A^{**}}(B^{**})A\subseteq B$, $\pi_r^{***}(b^{\prime\prime},e_{\alpha})\in B$. Consequently we have  $\pi_r^{***}(b^{\prime\prime},e_{\alpha}) \stackrel{w} {\rightarrow}\pi_r^{***}(b^{\prime\prime},e^{\prime\prime})=b^{\prime\prime}$ in $B$. It follows that $b^{\prime\prime}\in B$, since $B$ is a weakly complete.
\end{enumerate}
\end{proof}
\vspace{0.5cm}

\noindent A functional $a^\prime$ in $A^*$ is said to be $wap$ (weakly almost
 periodic) on $A$ if the mapping $a\rightarrow a^\prime a$ from $A$ into
 $A^{*}$ is weakly compact. The preceding definition is equivalent to the following condition, see [6, 14, 18].\\
 For any two net $(a_{\alpha})_{\alpha}$ and $(b_{\beta})_{\beta}$
 in $\{a\in A:~\parallel a\parallel\leq 1\}$, we have\\
$$\\lim_{\alpha}\\lim_{\beta}\langle  a^\prime,a_{\alpha}b_{\beta}\rangle  =\\lim_{\beta}\\lim_{\alpha}\langle  a^\prime,a_{\alpha}b_{\beta}\rangle  ,$$
whenever both iterated limits exist. The collection of all $wap$
functionals on $A$ is denoted by $wap(A)$. Also we have
$a^{\prime}\in wap(A)$ if and only if $\langle  a^{\prime\prime}b^{\prime\prime},a^\prime\rangle  =\langle  a^{\prime\prime}ob^{\prime\prime},a^\prime\rangle  $ for every $a^{\prime\prime},~b^{\prime\prime} \in
A^{**}$.\\\\

\noindent{\it{\bf Definition 2-7.}} Let $B$ be a Banach left $A-module$. Then, $b^\prime\in B^*$ is said to be left weakly almost periodic functional if the set $\{\pi^*_\ell(b^\prime,a):~a\in A,~\parallel a\parallel\leq 1\}$ is relatively weakly compact. We denote by $wap_\ell(B)$ the closed subspace of $B^*$ consisting of all the left weakly almost periodic functionals in $B^*$.\\
The definition of the right weakly almost periodic functional ($=wap_r(B)$) is the same.\\
By  [18], the definition of $wap_\ell(B)$ is equivalent to the following $$\langle  \pi_\ell^{***}(a^{\prime\prime},b^{\prime\prime}),b^\prime\rangle  =
\langle  \pi_\ell^{t***t}(a^{\prime\prime},b^{\prime\prime}),b^\prime\rangle  $$
for all $a^{\prime\prime}\in A^{**}$ and $b^{\prime\prime}\in B^{**}$.
Thus, we can write \\
$$wap_\ell(B)=\{ b^\prime\in B^*:~\langle  \pi_\ell^{***}(a^{\prime\prime},b^{\prime\prime}),b^\prime\rangle  =
\langle  \pi_\ell^{t***t}(a^{\prime\prime},b^{\prime\prime}),b^\prime\rangle  ~~$$$$for~~all~~a^{\prime\prime}\in A^{**},~b^{\prime\prime}\in B^{**}\}.$$\\
\noindent{\it{\bf Theorem 2-8.}} Suppose that $B$ is a Banach left $A-module$. Consider the following statements.
\begin{enumerate}

\item $B^*A\subseteq wap_\ell(B)$.
\item $AA^{**}\subseteq Z_{B^{**}}(A^{**})$.
\item $AA^{**}\subseteq AZ_{B^{**}}((A^{*}A)^*)$.
\end{enumerate}
Then, we have $(1)\Leftrightarrow (2)\Leftarrow (3)$.

\begin{proof} $(1)\Rightarrow (2)$\\
 Let $(b_{\alpha}^{\prime\prime})_{\alpha}\subseteq B^{**}$ such that  $b^{\prime\prime}_{\alpha} \stackrel{w^*} {\rightarrow}b^{\prime\prime}$. Then for all $a\in A$ and $a^{\prime\prime}\in A^{**}$, we have
$$\langle  \pi_\ell^{***}(aa^{\prime\prime},b^{\prime\prime}_{\alpha}),b^\prime\rangle  =
\langle  aa^{\prime\prime},\pi_\ell^{**}(b^{\prime\prime}_{\alpha},b^\prime)\rangle  =
\langle  a^{\prime\prime},\pi_\ell^{**}(b^{\prime\prime}_{\alpha},b^\prime)a\rangle  $$
$$=\langle  a^{\prime\prime},\pi_\ell^{**}(b^{\prime\prime}_{\alpha},b^\prime a)\rangle
=\langle  \pi_\ell^{***}(a^{\prime\prime},b^{\prime\prime}_{\alpha}),b^\prime a\rangle
\rightarrow\langle  \pi_\ell^{***}(a^{\prime\prime},b^{\prime\prime}),b^\prime a\rangle  $$
$$=\langle  \pi_\ell^{***}(aa^{\prime\prime},b^{\prime\prime}),b^\prime \rangle  .$$
Hence $aa^{\prime\prime}\in Z_{B^{**}}(A^{**})$.\\
$(2)\Rightarrow (1)$\\
Let $a\in A$ and $b^\prime\in B^*$. Then
$$\langle  \pi_\ell^{***}(a^{\prime\prime},b^{\prime\prime}),b^\prime a\rangle  =
\langle  a\pi_\ell^{***}(a^{\prime\prime},b^{\prime\prime}),b^\prime \rangle  =
\langle  \pi_\ell^{***}(aa^{\prime\prime},b^{\prime\prime}),b^\prime \rangle  $$
$$=\langle  \pi_\ell^{t***t}(aa^{\prime\prime},b^{\prime\prime}),b^\prime \rangle
=\langle  \pi_\ell^{t***t}(a^{\prime\prime},b^{\prime\prime}),b^\prime a\rangle  .$$
It follow that $b^\prime a\in wap_\ell(B)$.\\
$(3)\Rightarrow (2)$\\
Since $AZ_{B^{**}}((A^{*}A)^*)\subseteq Z_{B^{**}}(A^{**})$,  proof is hold.
\end{proof}

\vspace{0.5cm}

\noindent In the preceding theorem, if we take $B=A$, then we obtain Theorem 3.6 from [14] and  the same as preceding theorem, we can claim the following assertions:\\
If $B$ is a Banach right $A-module$, then for the following statements we have\\
 $(1)\Leftrightarrow (2)\Leftarrow (3)$.
 \begin{enumerate}

 \item $AB^*\subseteq wap_r(B)$.
\item $A^{**}A\subseteq Z_{B^{**}}(A^{**})$.
\item $A^{**}A\subseteq Z_{B^{**}}((A^{*}A)^*)A$.
\end{enumerate}
The proof of the this assertion is  similar to  proof of Theorem 2-8.\\\\

\noindent{\it{\bf Corollary 2-9.}} Suppose that $B$ is a  Banach $A-bimodule$. Then if $A$ is a left [resp. right] ideal in $A^{**}$, then $B^*A\subseteq wap_\ell(B)$ [resp. $AB^*\subseteq wap_r(B)$].\\\\

\noindent{\it{\bf Example 2-10.}} Suppose that $1\leq p\leq\infty$ and $q$ is conjugate of $p$. We know that if $G$ is compact, then $L^1(G)$ is a two-sided ideal in its second dual of it. By preceding Theorem we have  $L^q(G)*L^1(G)\subseteq wap_\ell(L^p(G))$ and $L^1(G)*L^q(G)\subseteq wap_r(L^p(G))$.\\
Also if $G$ is finite, then $L^q(G)\subseteq wap_\ell(L^p(G))\cap wap_r(L^p(G))$. Hence we conclude that\\
$$Z_{L^1(G)^{**}}(L^p(G)^{**})=L^p(G)~~~~ and~~Z_{L^p(G)^{**}}(L^1(G)^{**})=L^1(G).$$\\
\noindent{\it{\bf Theorem 2-11.}} We have the following assertions.
\begin{enumerate}
\item Suppose that $B$ is a Banach left $A-module$ and $b^\prime\in B^*$. Then $b^\prime\in wap_\ell(B)$ if and only if the adjoint of  the mapping $\pi_\ell^*(b^\prime,~):A\rightarrow B^*$ is $weak^*-to-weak$ continuous.
\item Suppose that $B$ is a Banach right $A-module$ and $b^\prime\in B^*$. Then $b^\prime\in wap_r(B)$ if and only if the adjoint of the mapping $\pi_r^*(b^\prime,~):B\rightarrow A^*$ is $weak^*-to-weak$ continuous.
\end{enumerate}
\begin{proof}
\begin{enumerate}
 \item Assume that $b^\prime\in wap_\ell(B)$ and $\pi_\ell^*(b^\prime,~)^*:B^{**}\rightarrow A^*$ is the adjoint of $\pi_\ell^*(b^\prime,~)$. Then for every $b^{\prime\prime}\in B^{**}$ and $a\in A$, we have
$$\langle  \pi_\ell^*(b^\prime,~)^*b^{\prime\prime},a\rangle  =\langle  b^{\prime\prime},\pi_\ell^*(b^\prime,a)\rangle  .$$
Suppose $(b_{\alpha}^{\prime\prime})_{\alpha}\subseteq B^{**}$ such that  $b^{\prime\prime}_{\alpha} \stackrel{w^*} {\rightarrow}b^{\prime\prime}$ and $a^{\prime\prime}\in A^{**}$ and $(a_{\beta})_{\beta}\subseteq A$ such that
 $a_{\beta}\stackrel{w^*} {\rightarrow}a^{\prime\prime}$.
 By easy calculation, for all $y^{\prime\prime}\in B^{**}$ and $y^\prime \in B^*$,  we have
  $$\langle  \pi_\ell^*(y^\prime,~)^*,y^{\prime\prime}\rangle  = \pi_\ell^{**}(y^{\prime\prime},y^\prime).$$
 Since $b^\prime\in wap_\ell(B)$,
$$\langle  \pi_\ell^{***}(a^{\prime\prime},b_{\alpha}^{\prime\prime}),b^{\prime}\rangle  \rightarrow \langle  \pi_\ell^{***}(a^{\prime\prime},b^{\prime\prime}),b^{\prime}\rangle  .$$
Then we have the following statements
$$\lim_\alpha\langle  a^{\prime\prime},\pi_\ell^*(b^\prime,~)^*b_{\alpha}^{\prime\prime}\rangle  =
\lim_\alpha\langle  a^{\prime\prime},\pi_\ell^{**}(b_{\alpha}^{\prime\prime},b^\prime)\rangle  $$
$$=
\lim_\alpha\langle  \pi_\ell^{***}(a^{\prime\prime},b_{\alpha}^{\prime\prime}),b^\prime\rangle  =
\langle  \pi_\ell^{***}(a^{\prime\prime},b^{\prime\prime}),b^\prime\rangle  $$
$$=\langle  a^{\prime\prime},\pi_\ell^*(b^\prime,~)^*b^{\prime\prime}\rangle  .$$
It follow that the adjoint of the mapping $\pi_\ell^*(b^\prime,~):A\rightarrow B^*$ is $weak^*-to-weak$ continuous.\\
Conversely, let the adjoint of the mapping $\pi_\ell^*(b^\prime,~):A\rightarrow B^*$ is $weak^*-to-weak$ continuous.
Suppose $(b_{\alpha}^{\prime\prime})_{\alpha}\subseteq B^{**}$ such that  $b^{\prime\prime}_{\alpha} \stackrel{w^*} {\rightarrow}b^{\prime\prime}$ and $b^\prime \in B^*$. Then for every $a^{\prime\prime}\in A^{**}$, we have
$$\lim_\alpha\langle  \pi_\ell^{***}(a^{\prime\prime},b_{\alpha}^{\prime\prime}),b^\prime\rangle  =
\lim_\alpha\langle  a^{\prime\prime},\pi_\ell^{**}(b_{\alpha}^{\prime\prime},b^\prime)\rangle  $$$$=
\lim_\alpha\langle  a^{\prime\prime},\pi_\ell^*(b^\prime,~)^*b_{\alpha}^{\prime\prime}\rangle  =
\langle  a^{\prime\prime},\pi_\ell^*(b^\prime,~)^*b^{\prime\prime}\rangle  =
\langle  \pi_\ell^{***}(a^{\prime\prime},b^{\prime\prime}),b^\prime\rangle  .$$
It follow that $b^\prime\in wap_\ell(B)$.\\
\item proof is similar to (1).
\end{enumerate}
\end{proof}
\vspace{0.5cm}

\noindent {\it{\bf Corollary 2-12.}} Let $A$ be a Banach algebra.
 Assume that $a^\prime\in A^*$ and $T_{a^\prime}$ is the linear operator from
$A$ into $A^*$ defined by $T_{a^\prime} a=a^\prime a$. Then, $a^\prime\in wap(A)$ if and
only if the adjoint of $T_{a^\prime}$ is
$weak^*-to-weak$ continuous. So $A$ is Arens regular if and only if the adjoint of the mapping  $T_{a^\prime} a=a^\prime a$ is $weak^*-to-weak$ continuous for every $a^\prime\in A^*$. \\\\\\

\begin{center}
\textbf{{ 3. $Lw^*w$-property and $Rw^*w$-property}}\\
\end{center}
\vspace{0.5cm}

\noindent In this section,  we introduce the new definition as $Left-weak^*-to-weak$ property and $Right-weak^*-to-weak$ property for Banach algebra $A$ and make some relations between these concepts and topological centers of module actions. As some conclusion, for locally compact group $G$, if $a\in M(G)$ has $Lw^*w-$ property [resp. $Rw^*w-$ property], then we have $L^1(G)^{**}*a\neq L^1(G)^{**}$ [resp. $a*L^1(G)^{**}\neq L^1(G)^{**}$], and also we have
 $Z_{L^1(G)^{**}}{(M(G)^{**})}\neq {M(G)^{**}}$ and $Z_{M(G)^{**}}{(L^1(G)^{**})}= {L^1(G)}$.  For a finite group $G$, we have $Z_{M(G)^{**}}{(L^1(G)^{**})}= {L^1(G)^{**}}$ and
$Z_{L^1(G)^{**}}{(M(G)^{**})}= {M(G)^{**}}$.\\\\

\noindent{\it{\bf Definition 3-1.}} Let $B$ be a Banach left $A-module$. We say that $a\in A$ has $Left-weak^*-to-weak$ property ($=Lw^*w-$ property) with respect to $B$, if for all $(b^\prime_{\alpha})_{\alpha}\subseteq B^*$, $ab_\alpha^\prime\stackrel{w^*} {\rightarrow}0$ implies $ab_\alpha^\prime\stackrel{w} {\rightarrow}0$. If every $a\in A$ has $Lw^*w-$ property with respect to $B$, then we say that $A$ has $Lw^*w-$ property with respect to $B$. The definition of the $Right-weak^*-to-weak$ property ($=Rw^*w-$ property) is the same.\\
We say that $a\in A$ has $weak^*-to-weak$  property ($=w^*w-$ property) with respect to $B$ if it has $Lw^*w-$ property and $Rw^*w-$ property with respect to $B$.\\
If $a\in A$ has $Lw^*w-$ property with respect to itself, then we say that $a\in A$ has $Lw^*w-$ property.\\\\
For preceding definition, we have  some examples and remarks as follows.\\
a) If $B$ is a Banach $A$-bimodule and reflexive, then $A$ has $w^*w-$property with respect to $B$. Then we have the following statenents for group algebras.\\
i) $L^1(G)$, $M(G)$ and $A(G)$ have $w^*w-$property when $G$ is finite.\\
ii) Let $G$ be locally compact group.  $L^1(G)$ [resp. $M(G)$] has $w^*w-$property [resp. $Lw^*w-$ property ] with respect to $L^p(G)$ whenever $p> 1$.\\
b) Suppose that $B$ is a Banach left $A-module$ and $e$ is left unit element of $A$ such that $eb=b$ for all $b\in B$. If $e$ has $Lw^*w-$ property, then $B$ is reflexive.\\
c) If $S$ is a compact semigroup, then $C^+(S)=\{f\in C(S):~f>0\}$ has $w^*w-$property.\\\\\\
\noindent{\it{\bf Theorem 3-2.}} Suppose that $B$ is a Banach  $A-bimodule$. Then we have the following assertions.
 \begin{enumerate}
\item If $A^{**}=a_0A^{**}$ [resp.  $A^{**}=A^{**}a_0$] for some $a_0\in A$ and $a_0$ has $Rw^*w-$ property [resp. $Lw^*w-$ property], then $Z_{B^{**}}(A^{**})=A^{**}$.
\item If $B^{**}=a_0B^{**}$ [resp.  $B^{**}=B^{**}a_0$] for some $a_0\in A$ and $a_0$ has $Rw^*w-$ property [resp. $Lw^*w-$ property] with respect to $B$, then $Z_{A^{**}}(B^{**})=B^{**}$.
\end{enumerate}
\begin{proof}
\begin{enumerate}
\item Suppose that $A^{**}=a_0A^{**}$  for some $a_0\in A$ and $a_0$ has $Rw^*w-$ property. Let   $(b_{\alpha}^{\prime\prime})_{\alpha}\subseteq B^{**}$ such that  $b^{\prime\prime}_{\alpha} \stackrel{w^*} {\rightarrow}b^{\prime\prime}$. Then for all $a\in A$ and $b^\prime\in B^{*}$, we have
$$\langle  \pi_\ell^{**}(b_\alpha^{\prime\prime},b^{\prime}), a\rangle  =\langle  b_\alpha^{\prime\prime},\pi_\ell^{*}(b^{\prime}, a)\rangle  \rightarrow
\langle  b^{\prime\prime},\pi_\ell^{*}(b^{\prime}, a)\rangle  =\langle  \pi_\ell^{**}(b^{\prime\prime},b^{\prime}), a\rangle  ,$$
it follow that $\pi_\ell^{**}(b_\alpha^{\prime\prime},b^{\prime})\stackrel{w^*} {\rightarrow}\pi_\ell^{**}(b^{\prime\prime},b^{\prime})$. Also we can write  $\pi_\ell^{**}(b_\alpha^{\prime\prime},b^{\prime})a_0\stackrel{w^*} {\rightarrow}\pi_\ell^{**}(b^{\prime\prime},b^{\prime})a_0$. Since $a_0$ has $Rw^*w-$ property, $\pi_\ell^{**}(b_\alpha^{\prime\prime},b^{\prime})a_0\stackrel{w} {\rightarrow}\pi_\ell^{**}(b^{\prime\prime},b^{\prime})a_0$. Now let $a^{\prime\prime}\in A^{**}$. Then there is
$x^{\prime\prime}\in A^{**}$ such that  $a^{\prime\prime}=a_0x^{\prime\prime}$ consequently we have
$$\langle  \pi_\ell^{***}(a^{\prime\prime},b^{\prime\prime}_{\alpha}),b^\prime\rangle  =
\langle  a^{\prime\prime},\pi_\ell^{**}(b^{\prime\prime}_{\alpha},b^\prime)\rangle  =
\langle  x^{\prime\prime},\pi_\ell^{**}(b^{\prime\prime}_{\alpha},b^\prime)a_0\rangle  $$
$$\rightarrow\langle  x^{\prime\prime},\pi_\ell^{**}(b^{\prime\prime},b^\prime)a_0\rangle  =
\langle  \pi_\ell^{***}(a^{\prime\prime},b^{\prime\prime}_{\alpha}),b^\prime\rangle  .$$
We conclude that $a^{\prime\prime}\in Z_{B^{**}}(A^{**})$. Proof of the next part is the same as the preceding proof.\\
\item Let $B^{**}=a_0B^{**}$  for some $a_0\in A$ and $a_0$ has $Rw^*w-$ property  with respect to $B$. Assume that
$(a_{\alpha}^{\prime\prime})_{\alpha}\subseteq A^{**}$ such that  $a^{\prime\prime}_{\alpha} \stackrel{w^*} {\rightarrow}a^{\prime\prime}$. Then for all $b\in B$, we have
$$\langle  \pi_r^{**}(a_\alpha^{\prime\prime},b^{\prime}), b\rangle  =\langle  a_\alpha^{\prime\prime},\pi_r^{*}(b^{\prime}, b)\rangle
\rightarrow\langle  a^{\prime\prime},\pi_r^{*}(b^{\prime}, b)\rangle  =\langle  \pi_r^{**}(a^{\prime\prime},b^{\prime}), b\rangle  .$$
We conclude that $\pi_r^{**}(a_\alpha^{\prime\prime},b^{\prime})\stackrel{w^*} {\rightarrow}
\pi_r^{**}(a^{\prime\prime},b^{\prime})$ then we have $\pi_r^{**}(a_\alpha^{\prime\prime},b^{\prime})a_0\stackrel{w^*} {\rightarrow}
\pi_r^{**}(a^{\prime\prime},b^{\prime})a_0$. Since $a_0$ has $Rw^*w-$ property with respect to $B$,
$\pi_r^{**}(a_\alpha^{\prime\prime},b^{\prime})a_0\stackrel{w} {\rightarrow}
\pi_r^{**}(a^{\prime\prime},b^{\prime})a_0$.\\
Now let $b^{\prime\prime}\in B^{**}$. Then there is $x^{\prime\prime}\in B^{**}$ such that $b^{\prime\prime}=
a_0x^{\prime\prime}$. Hence, we have
$$\langle  \pi_r^{***}(b^{\prime\prime},a_\alpha^{\prime\prime}), b^\prime\rangle  =\langle  b^{\prime\prime},\pi_r^{**}(a^{\prime\prime}_\alpha, b^\prime)\rangle  =\langle  a_0x^{\prime\prime},\pi_r^{**}(a^{\prime\prime}_\alpha, b^\prime)\rangle  $$
$$=\langle  x^{\prime\prime},\pi_r^{**}(a^{\prime\prime}_\alpha, b^\prime)a_0\rangle  \rightarrow \langle  x^{\prime\prime},\pi_r^{**}(a^{\prime\prime}, b^\prime)a_0\rangle  =\langle  b^{\prime\prime},\pi_r^{**}(a^{\prime\prime}, b^\prime)\rangle  $$
$$=\langle  \pi_r^{***}(b^{\prime\prime},a^{\prime\prime}), b^\prime\rangle  .$$
It follow that $b^{\prime\prime}\in Z_{A^{**}}(B^{**})$. The next part is similar to the preceding proof.
\end{enumerate}
\end{proof}
\vspace{0.5cm}

\noindent{\it{\bf Example 3-3.}}
\begin{enumerate}

\item By using Theorem 2-3, for locally compact group $G$, if $a\in M(G)$ has $Lw^*w-$ property [resp. $Rw^*w-$ property], then we have $L^1(G)^{**}*a\neq L^1(G)^{**}$ [resp. $a*L^1(G)^{**}\neq L^1(G)^{**}$].

\item  If $G$ is finite, then by Theorem 2-3, we have $Z_{M(G)^{**}}{(L^1(G)^{**})}= {L^1(G)^{**}}$ and
$Z_{L^1(G)^{**}}{(M(G)^{**})}= {M(G)^{**}}$.\\\\
\end{enumerate}

\noindent{\it{\bf Theorem 3-4.}} Suppose that $B$ is a Banach  $A-bimodule$ and $A$ has a $BAI$. Then we have the following assertions.
 \begin{enumerate}
\item  If $B^*$ factors on the left [resp. right] with respect to $A$ and $A$ has $Rw^*w-$ property [resp. $Lw^*w-$ property], then $Z_{B^{**}}(A^{**})=A^{**}$.
\item  If $B^*$ factors on the left [resp. right] with respect to $A$ and $A$ has $Rw^*w-$ property [resp. $Lw^*w-$ property] with respect $B$, then $Z_{A^{**}}(B^{**})=B^{**}$.
\end{enumerate}
\begin{proof}
 \begin{enumerate}
\item   Assume that $B^*$ factors on the left and $A$ has $Rw^*w-$ property. Let   $(b_{\alpha}^{\prime\prime})_{\alpha}\subseteq B^{**}$ such that  $b^{\prime\prime}_{\alpha} \stackrel{w^*} {\rightarrow}b^{\prime\prime}$. Since $B^*A=B^*$, for all $b^\prime\in B^*$ there are $x\in A$ and $y^\prime\in B^*$ such that $b^\prime=y^\prime x$. Then for all $a\in A$, we have
$$\langle  \pi_\ell^{**}(b_\alpha^{\prime\prime},y^{\prime})x,a\rangle  =\langle  b_\alpha^{\prime\prime},\pi_\ell^{*}(y^{\prime},a)x\rangle
=\langle  \pi_\ell^{**}(b_\alpha^{\prime\prime},b^{\prime}),a\rangle  $$$$=\langle  b_\alpha^{\prime\prime},\pi_\ell^{*}(b^{\prime},a)\rangle
\rightarrow\langle  b^{\prime\prime},\pi_\ell^{*}(b^{\prime},a)\rangle
=\langle  \pi_\ell^{**}(b^{\prime\prime},y^{\prime})x,a\rangle  .$$
Thus, we conclude that $\pi_\ell^{**}(b_\alpha^{\prime\prime},y^{\prime})x\stackrel{w^*} {\rightarrow}
\pi_\ell^{**}(b^{\prime\prime},y^{\prime})x$.
Since $A$ has $Rw^*w-$ property, $\pi_\ell^{**}(b_\alpha^{\prime\prime},y^{\prime})x\stackrel{w} {\rightarrow}
\langle  \pi_\ell^{**}(b^{\prime\prime},y^{\prime})x$.
Now let $a^{\prime\prime}\in A^{**}$. Then
$$\langle  \pi_\ell^{***}(a^{\prime\prime},b^{\prime\prime}_{\alpha}),b^\prime\rangle  =
\langle  a^{\prime\prime},\pi_\ell^{**}(b^{\prime\prime}_{\alpha},b^\prime)\rangle  =
\langle  a^{\prime\prime},\pi_\ell^{**}(b^{\prime\prime}_{\alpha},y^\prime)x\rangle  $$$$\rightarrow
\langle  a^{\prime\prime},\pi_\ell^{**}(b^{\prime\prime},y^\prime)x\rangle  =
\langle  \pi_\ell^{***}(a^{\prime\prime},b^{\prime\prime}),b^\prime\rangle  .$$
It follow that $a^{\prime\prime}\in Z_{B^{**}}(A^{**})$.\\
If $B^*$ factors on the right with respect to $A$, and assume that  $A$ has $Lw^*w-$ property, then proof is the same as preceding proof.\\
\item  Let $B^*$ factors on the left  with respect to $A$ and $A$ has $Rw^*w-$ property  with respect to $B$.
Assume that
$(a_{\alpha}^{\prime\prime})_{\alpha}\subseteq A^{**}$ such that  $a^{\prime\prime}_{\alpha} \stackrel{w^*} {\rightarrow}a^{\prime\prime}$. Since $B^*A=B^*$, for all $b^\prime\in B^*$ there are $x\in A$ and $y^\prime\in B^*$ such that $b^\prime=y^\prime x$. Then for all $b\in B$, we have
$$\langle  \pi_r^{**}(a_\alpha^{\prime\prime},y^{\prime})x, b\rangle  =\langle  \pi_r^{**}(a_\alpha^{\prime\prime},b^{\prime}), b\rangle  =
\langle  a_\alpha^{\prime\prime},\pi_r^{*}(b^{\prime}, b)\rangle  $$$$=\langle  a^{\prime\prime},\pi_r^{*}(b^{\prime}, b)\rangle  =\langle  \pi_r^{**}(a^{\prime\prime},y^{\prime})x, b\rangle  .$$
Consequently  $\pi_r^{**}(a_\alpha^{\prime\prime},y^{\prime})x\stackrel{w^*} {\rightarrow}\pi_r^{**}(a^{\prime\prime},y^{\prime})x$. Since $A$ has $Rw^*w-$ property  with respect to $B$,
$\pi_r^{**}(a_\alpha^{\prime\prime},y^{\prime})x\stackrel{w} {\rightarrow}\pi_r^{**}(a^{\prime\prime},y^{\prime})x$.
It follow that for all $b^{\prime\prime}\in B^{**}$, we have
$$\langle  \pi_r^{***}(b^{\prime\prime},a_\alpha^{\prime\prime}), b^\prime\rangle  =
\langle  b^{\prime\prime},\pi_r^{**}(a_\alpha^{\prime\prime}, y^\prime)x\rangle  \rightarrow
\langle  b^{\prime\prime},\pi_r^{**}(a^{\prime\prime}, y^\prime)x\rangle  $$
$$=\langle  \pi_r^{***}(b^{\prime\prime},a^{\prime\prime}), b^\prime\rangle  .$$
Thus we conclude that $b^{\prime\prime}\in Z_{A^{**}}(B^{**})$.\\
The proof of the next assertions is the same as the preceding proof.
\end{enumerate}
\end{proof}
\vspace{0.5cm}

\noindent{\it{\bf Theorem 3-5.}} Suppose that $B$ is a Banach  $A-bimodule$. Then we have the following assertions.
\begin{enumerate}
\item  If $a_0\in A$ has $Rw^*w-$ property with respect to $B$, then $a_0A^{**}\subseteq Z_{B^{**}}(A^{**})$ and $a_0B^*\subseteq wap_\ell(B)$.
\item  If $a_0\in A$ has $Lw^*w-$ property with respect to $B$, then $A^{**}a_0\subseteq Z_{B^{**}}(A^{**})$ and $B^*a_0\subseteq wap_\ell(B)$.
\item  If $a_0\in A$ has $Rw^*w-$ property with respect to $B$, then $a_0B^{**}\subseteq Z_{A^{**}}(B^{**})$ and $B^*a_0\subseteq wap_r(B)$.
\item  If $a_0\in A$ has $Lw^*w-$ property with respect to $B$, then $B^{**}a_0\subseteq Z_{A^{**}}(B^{**})$ and $a_0B^*\subseteq wap_r(B)$.
\end{enumerate}
\begin{proof}
\begin{enumerate}

\item Let   $(b_{\alpha}^{\prime\prime})_{\alpha}\subseteq B^{**}$ such that  $b^{\prime\prime}_{\alpha} \stackrel{w^*} {\rightarrow}b^{\prime\prime}$. Then for all $a\in A$ and  $b^\prime\in B^*$, we have
$$\langle  \pi_\ell^{**}(b_\alpha^{\prime\prime},b^{\prime})a_0,a\rangle  =\langle  \pi_\ell^{**}(b_\alpha^{\prime\prime},b^{\prime}),a_0a\rangle
=\langle  b_\alpha^{\prime\prime},\pi_\ell^{*}(b^{\prime},a_0a)\rangle  $$
$$\rightarrow \langle  b^{\prime\prime},\pi_\ell^{*}(b^{\prime},a_0a)\rangle  =\langle  \pi_\ell^{**}(b^{\prime\prime},b^{\prime})a_0,a\rangle  .$$
It follow that $\pi_\ell^{**}(b_\alpha^{\prime\prime},b^{\prime})a_0\stackrel{w^*} {\rightarrow}
\pi_\ell^{**}(b^{\prime\prime},b^{\prime})a_0$. Since $a_0$ has $Rw^*w-$ property with respect to $B$,
$\pi_\ell^{**}(b_\alpha^{\prime\prime},b^{\prime})a_0\stackrel{w} {\rightarrow}
\pi_\ell^{**}(b^{\prime\prime},b^{\prime})a_0$.\\
We conclude that $a_0a^{\prime\prime}\in  Z_{B^{**}}(A^{**})$ so that   $a_0 A^{**}\in  Z_{B^{**}}(A^{**})$.
Since $\pi_\ell^{**}(b^{\prime\prime},b^{\prime})a_0=\pi_\ell^{**}(b^{\prime\prime},b^{\prime}a_0)$, $a_0B^*\subseteq wap_\ell(B)$.\\
\item proof is similar to (1).
\item Assume that
$(a_{\alpha}^{\prime\prime})_{\alpha}\subseteq A^{**}$ such that  $a^{\prime\prime}_{\alpha} \stackrel{w^*} {\rightarrow}a^{\prime\prime}$. Let $b\in B$ and $b^\prime\in B^*$. Then we have
$$\langle  \pi_r^{**}(a_\alpha^{\prime\prime},b^{\prime})a_0, b\rangle  =\langle  \pi_r^{**}(a_\alpha^{\prime\prime},b^{\prime}),a_0 b\rangle  =
\langle  a_\alpha^{\prime\prime},\pi_r^{*}(b^{\prime},a_0 b)\rangle  $$
$$\rightarrow \langle  a^{\prime\prime},\pi_r^{*}(b^{\prime},a_0 b)\rangle  =\langle  \pi_r^{**}(a^{\prime\prime},b^{\prime})a_0, b\rangle  .$$
Thus we conclude $\pi_r^{**}(a_\alpha^{\prime\prime},b^{\prime})a_0\stackrel{w^*} {\rightarrow}\pi_r^{**}(a^{\prime\prime},b^{\prime})a_0$.
Since $a_0$ has $Rw^*w-$ property with respect to $B$, $\pi_r^{**}(a_\alpha^{\prime\prime},b^{\prime})a_0\stackrel{w} {\rightarrow}\pi_r^{**}(a^{\prime\prime},b^{\prime})a_0$.
If $b^{\prime\prime}\in B^{**}$, then we have
$$\langle  \pi_r^{***}(a_0b^{\prime\prime},a_\alpha^{\prime\prime}), b^\prime\rangle  =\langle  a_0b^{\prime\prime},\pi_r^{***}(a_\alpha^{\prime\prime}, b^\prime)\rangle  =
\langle  b^{\prime\prime},\pi_r^{**}(a_\alpha^{\prime\prime}, b^\prime)a_0\rangle  $$
$$=\langle  b^{\prime\prime},\pi_r^{**}(a_\alpha^{\prime\prime}, b^\prime)a_0\rangle  =\langle  \pi_r^{***}(a_0b^{\prime\prime},a^{\prime\prime}), b^\prime\rangle  .$$
It follow that $a_0b^{\prime\prime}\in Z_{A^{**}}(B^{**})$. Consequently we have  $a_0B^{**}\in Z_{A^{**}}(B^{**})$.\\
The proof of the next assertion is clear.\\
\item Proof is similar to (3).
\end{enumerate}
\end{proof}
\vspace{0.5cm}

\noindent{\it{\bf Theorem 3-6.}} Let $B$ be a Banach  $A-bimodule$. Then we have the following assertions.
\begin{enumerate}
\item  Suppose
$$\lim_{\alpha}\lim_{\beta}\langle  b^\prime_{\beta},b_{\alpha}\rangle  =\lim_{\beta}\lim_{\alpha}
\langle  b^\prime_{\beta},b_{\alpha}\rangle  ,$$ for every $(b_{\alpha})_{\alpha}\subseteq B$ and $(b^\prime_{\beta})_{\beta}\subseteq B^*$. Then $A$ has $Lw^*w-$ property and $Rw^*w-$ property with respect to $B$.\\
\item  If for some $a\in A$, $$\lim_{\alpha}\lim_{\beta}\langle  ab^\prime_{\beta},b_{\alpha}\rangle  =\lim_{\beta}\lim_{\alpha}
\langle  ab^\prime_{\beta},b_{\alpha}\rangle  ,$$
 for every $(b_{\alpha})_{\alpha}\subseteq B$ and $(b^\prime_{\beta})_{\beta}\subseteq B^*$, then $a$ has $Rw^*w-$ property with respect to $B$. Also if for some $a\in A$, $$\lim_{\alpha}\lim_{\beta}\langle  b^\prime_{\beta}a,b_{\alpha}\rangle  =\lim_{\beta}\lim_{\alpha}
\langle  b^\prime_{\beta}a,b_{\alpha}\rangle  ,$$
for every $(b_{\alpha})_{\alpha}\subseteq B$ and $(b^\prime_{\beta})_{\beta}\subseteq B^*$, then $a$ has $Lw^*w-$ property with respect to $B$.
\end{enumerate}
\begin{proof}
\begin{enumerate}
 \item  Assume that $a\in A$ such that $ab_{\beta}^\prime\stackrel{w^*} {\rightarrow}0$ where
$(b^\prime_{\beta})_{\beta}\subseteq B^*$. Let $b^{\prime\prime}\in B^{**}$ and $(b_{\alpha})_{\alpha}\subseteq B$
such that $b_{\alpha}\stackrel{w^*} {\rightarrow}b^{\prime\prime}$. Then
$$\lim_{\beta}\langle  b^{\prime\prime},ab_{\beta}^\prime\rangle  =\lim_{\beta}\lim_{\alpha}\langle  b_{\alpha},ab_{\beta}^\prime\rangle  =
\lim_{\beta}\lim_{\alpha}\langle  ab_{\beta}^\prime ,b_{\alpha}\rangle  $$
$$=\lim_{\alpha}\lim_{\beta}\langle  ab_{\beta}^\prime ,b_{\alpha}\rangle  =0.$$
We conclude that $ab_{\beta}^\prime\stackrel{w} {\rightarrow}0$, so $A$ has $Lw^*w-$ property. It also easy that $A$ has  $Rw^*w-$ property.\\
\item  Proof is easy and is the same as (1).
\end{enumerate}\end{proof}
\vspace{0.5cm}

\noindent{\it{\bf Definition 3-7.}} Let $B$ be a Banach left $A-module$. We say that $B^*$ strong factors on the left [resp. right] if for all $(b^\prime_\alpha)_\alpha\subseteq B^*$ there are $(a_\alpha)_\alpha\subseteq A$ and $b^\prime\in B^*$ such that $b^\prime_\alpha=b^\prime a_\alpha$ [resp. $b^\prime_\alpha= a_\alpha b^\prime$] where $(a_\alpha)_\alpha$ has limit the $weak^*$ topology in $A^{**}$.\\
If $B^*$ strong factors on the left and right, then we say that $B^*$ strong factors on the both side.\\
It is clear that if $B^*$ strong factors on the left [resp. right], then $B^*$ factors on the left [resp. right].\\\\

\noindent{\it{\bf Theorem 3-8.}} Suppose that $B$ is a Banach  $A-bimodule$. Assume that $AB^*\subseteq wap_\ell B$. If
$B^*$ strong factors on the left [resp. right], then $A$ has $Lw^*w-$ property [resp. $Rw^*w-$ property ] with respect to $B$.
\begin{proof}  Let   $(b_{\alpha}^{\prime})_{\alpha}\subseteq B^{*}$ such that  $ab_\alpha^\prime\stackrel{w^*} {\rightarrow}0$. Since $B^*$ strong factors on the left, there are
$(a_\alpha)_\alpha\subseteq A$ and $b^\prime\in B^*$ such that $b^\prime_\alpha=b^\prime a_\alpha$. Let
$b^{\prime\prime}\in B^{**}$ and $(b_{\beta})_{\beta}\subseteq B$ such that
 $b_{\beta}\stackrel{w^*} {\rightarrow}b^{\prime\prime}$. Then we have
$$\lim_\alpha\langle  b^{\prime\prime},ab_\alpha^\prime\rangle  =\lim_\alpha\lim_{\beta}\langle  b_{\beta},ab_\alpha^\prime\rangle
=\lim_\alpha\lim_{\beta}\langle  ab_\alpha^\prime,b_{\beta}\rangle  $$
$$=\lim_\alpha\lim_{\beta}\langle  ab^\prime a_\alpha,b_{\beta}\rangle  =\lim_\alpha\lim_{\beta}\langle  ab^\prime, a_\alpha b_{\beta}\rangle  $$
$$=\lim_{\beta}\lim_\alpha\langle  ab^\prime, a_\alpha b_{\beta}\rangle  =\lim_{\beta}\lim_\alpha\langle  ab_\alpha^\prime,b_{\beta}\rangle  =0$$
It follow that $ab_\alpha^\prime\stackrel{w} {\rightarrow}0$.
\end{proof}

\vspace{0.5cm}
\noindent{\it{\bf Problems .}}
\begin{enumerate}
 \item Suppose that $B$ is a Banach  $A-bimodule$. If $B$ is left or right factors with respect to $A$, dose $A$ has $Lw^*w-$property or $Rw^*w-$property, respectively?
\item Suppose that $B$ is a Banach  $A-bimodule$. Let $A$ has $Lw^*w-$property with respect to $B$. Dose $Z_{B^{**}}(A^{**})=A^{**}$?\\\\\\
\end{enumerate}

\begin{center}
\textbf{{ 4. Factorization  properties and  topological centers of \\ left module actions}}
\end{center}
\vspace{0.5cm}

Let $B$ be a left Banach $A-module$ and  $e$ be a left  unit element of $A$. We say that $e$ is a left unit (resp. weakly left unit)  $A-module$ for $B$, if $\pi_\ell(e,b)=b$ (resp. $\langle  b^\prime , \pi_\ell(e,b)\rangle=\langle  b^\prime , b\rangle$ for all $b^\prime\in B^*$) where $b\in B$. The definition of right unit (resp. weakly right unit) $A-module$ is similar.\\
We say that a Banach $A-bimodule$ $B$ is a unital $A-module$, if $B$ has left and right unit $A-module$ that are equal then we say that $B$ is unital $A-module$.\\
Let $B$ be a left Banach $A-module$ and $(e_{\alpha})_{\alpha}\subseteq A$ be a LAI [resp. weakly left approximate identity(=WLAI)] for $A$. We say that $(e_{\alpha})_{\alpha}$ is left approximate identity  ($=LAI$)[ resp. weakly left approximate identity  (=$WLAI$)] for $B$, if for all $b\in B$,  $\pi_\ell (e_{\alpha},b) \stackrel{} {\rightarrow}
b$ ( resp. $\pi_\ell (e_{\alpha},b) \stackrel{w} {\rightarrow}
b$). The definition of the right approximate identity ($=RAI$)[ resp. weakly right approximate identity ($=WRAI$)] is similar.\\
 $(e_{\alpha})_{\alpha}\subseteq A$ is called a approximate identity  ($=AI$)[ resp. weakly approximate identity  ($WAI$)] for $B$, if $B$ has the same left and right approximate identity  [ resp. weakly left and right approximate identity ].\\
Let $(e_{\alpha})_{\alpha}\subseteq A$ be $weak^*$ left approximate identity for $A^{**}$. We say that $(e_{\alpha})_{\alpha}$ is $weak^*$ left approximate identity as $A^{**}-module$ ($=W^*LAI~as~A^{**}-module$) for $B^{**}$, if for all $b^{\prime\prime}\in B^{**}$, we have $\pi_\ell^{***} (e_{\alpha},b^{\prime\prime}) \stackrel{w^*} {\rightarrow}
b^{\prime\prime}$. The definition of the $weak^*$ right approximate identity ($=W^*RAI~$) is similar.\\
 $(e_{\alpha})_{\alpha}\subseteq A$ is called a $weak^*$ approximate identity ($=W^*AI~$) for $B^{**}$, if $B^{**}$ has the same $weak^*$ left and right approximate identity.\\
Let $B$ be a Banach $A-bimodule$. We say that $B$ is a left [resp. right] factors with respect to $A$, if $BA=B$ [resp. $AB=B$].\\\\

\noindent{\it{\bf Theorem 4-1.}} Let $B$ be a Banach left $A-module$ and $(e_{\alpha})_{\alpha}\subseteq A$ be a $LBAI$ for $B$. Then the following assertions hold.
\begin{enumerate}
\item For each $b^\prime\in B^*$, we have  $\pi_\ell^*(b^\prime,e_\alpha)\stackrel{w^*} {\rightarrow}b^\prime$.
\item $B^*$ factors on the left with respect to $A$ if and only if $B^{**}$ has a $W^*LBAI$ $(e_{\alpha})_{\alpha}\subseteq A$.
\item $B^{**}$ has a $W^*LBAI$ $(e_{\alpha})_{\alpha}\subseteq A$ if and only if $B^{**}$ has a left unit element $e^{\prime \prime}\in A^{**}$ such that $e_{\alpha}\stackrel{w^*} {\rightarrow}e^{\prime \prime}$.
\end{enumerate}
\begin{proof}
\begin{enumerate}
\item Let $b\in B$ and $b^\prime\in B^*$. Since $\mid \langle  b^\prime , \pi_\ell(e_\alpha ,b)\mid\leq \parallel b^\prime\parallel\parallel \pi_\ell(e_\alpha ,b)\parallel$, we have the following equality
$$\lim_\alpha \langle  \pi^*_\ell(b^\prime,e_\alpha),b\rangle=\lim_\alpha\langle  b^\prime,\pi_\ell(e_\alpha,b)\rangle= 0.$$
It follows that  $\pi^*_\ell(b^\prime,e_\alpha)\stackrel{w^*} {\rightarrow}0$.

\item  Let $B^*$ factors on the left with respect to $A$. Then for every $b^\prime\in B^*$, there are $x^\prime\in B^*$ and $a\in A$ such that $b^\prime=x^\prime a$. Then for every $b^{\prime\prime}\in B^{**}$, we have
$$\langle  \pi^{***}_\ell(e_\alpha,b^{\prime\prime}),b^\prime \rangle= \langle  e_\alpha,\pi^{**}_\ell(b^{\prime\prime},b^\prime )\rangle=
\langle  \pi^{**}_\ell(b^{\prime\prime},b^\prime ),e_\alpha\rangle$$$$=\langle  b^{\prime\prime},\pi^{*}_\ell(b^\prime ,e_\alpha)\rangle=
\langle  b^{\prime\prime},\pi^{*}_\ell(x^\prime a,e_\alpha)\rangle= \langle  b^{\prime\prime},\pi^{*}_\ell(x^\prime ,ae_\alpha)\rangle $$$$=
 \langle  \pi^{**}_\ell(b^{\prime\prime},x^\prime ),ae_\alpha\rangle \rightarrow   \langle  \pi^{**}_\ell(b^{\prime\prime},x^\prime ),a\rangle
 $$$$=\langle   b^{\prime\prime},b^\prime \rangle .$$
It follows that
$$\pi^{***}_\ell(e_\alpha,b^{\prime\prime})\stackrel{w^*} {\rightarrow}b^{\prime\prime},$$
 and so $B^{**}$ has $W^*LBAI$.\\
Conversely, let $b^\prime\in B^*$. Then for every $b^{\prime\prime}\in B^{**}$, we have
$$\langle  b^{\prime\prime},\pi^{*}_\ell(b^\prime ,e_\alpha)\rangle =\langle  \pi^{***}_\ell(e_\alpha,b^{\prime\prime}),b^\prime \rangle \rightarrow\langle   b^{\prime\prime},b^\prime \rangle .$$
It follows that
$$\pi^{*}_\ell(b^\prime ,e_\alpha)\stackrel{w} {\rightarrow} b^\prime ,$$
and so by Cohen factorization theorem, we are done.

\item Assume that  $B^{**}$ has a $W^*LBAI$ $(e_{\alpha})_{\alpha}\subseteq A$. Without loss generality, let $e^{\prime\prime}\in A^{**}$ be a left unit for $A^{**}$ with respect to the first Arens product such that $e_{\alpha}\stackrel{w^*} {\rightarrow}e^{\prime\prime}$. Then for each $b^\prime\in B^*$, we have

 $$\langle  \pi^{***}_\ell(e^{\prime\prime},b^{\prime\prime}),b^\prime \rangle =\langle  e^{\prime\prime},\pi^{**}_\ell(b^{\prime\prime},b^\prime ) \rangle $$$$=\lim_\alpha\langle  e_{\alpha},\pi^{**}_\ell(b^{\prime\prime},b^\prime) \rangle =\lim_\alpha\langle  \pi^{**}_\ell(b^{\prime\prime},b^\prime ),e_{\alpha}\rangle
 $$$$=\lim_\alpha\langle  b^{\prime\prime},\pi^{*}_\ell(b^\prime ,e_{\alpha})\rangle
=\lim_\alpha\langle  \pi^{****}_\ell(b^\prime ,e_{\alpha}),b^{\prime\prime}\rangle
$$$$=\lim_\alpha\langle  b^\prime ,\pi^{***}_\ell(e_{\alpha},b^{\prime\prime})\rangle
=\lim_\alpha\langle  \pi^{***}_\ell(e_{\alpha},b^{\prime\prime}),b^\prime \rangle $$$$= \langle   b^{\prime\prime},b^\prime \rangle .$$
Thus $e^{\prime\prime}\in A^{**}$ is a left unit for $B^{**}$.\\
Conversely, let  $e^{\prime\prime}\in A^{**}$ be a left unit for $B^{**}$ and assume that   $e_{\alpha}\stackrel{w^*} {\rightarrow}e^{\prime\prime}$ in $A^{**}$. Then for  every $b^{\prime\prime}\in B^{**}$ and $b^\prime\in B^*$, we have
$$\langle  \pi^{***}_\ell(e_\alpha,b^{\prime\prime}),b^\prime \rangle =\langle  e_\alpha,\pi^{**}_\ell(b^{\prime\prime},b^\prime )\rangle $$$$\rightarrow
\langle  e^{\prime\prime},\pi^{**}_\ell(b^{\prime\prime},b^\prime )\rangle =\langle  \pi^{***}_\ell(e^{\prime\prime},b^{\prime\prime}),b^\prime \rangle $$
$$= \langle   b^{\prime\prime},b^\prime \rangle .$$\\

It follows that $\pi^{***}_\ell(e_\alpha,b^{\prime\prime})\stackrel{w^*} {\rightarrow}b^{\prime\prime}.$\end{enumerate}\end{proof}

\vspace{0.5cm}
\noindent{\it{\bf Corollary 4-2.}} Let $B$ be a Banach left $A-module$ and $A$ has a $BLAI$. If $B^{**}$ has a $W^*LBAI$ , then
$$\{a^{\prime\prime}\in A^{**}:~Aa^{\prime\prime}\subseteq A\}\subseteq Z_{B^{**}}(A^{**}).$$
\begin{proof} By using the preceding theorem, since $B^{**}$ has $W^*LBAI$, $B^*$ factors on the left with respect to $A$. Suppose that $b^\prime\in B^*$. Then there are $x^\prime\in B^*$ and $a\in A$ such that $b^\prime=x^\prime a$. Assume that $a^{\prime\prime}\in A^{**}$ such that $Aa^{\prime\prime}\subseteq A$. Let $b^{\prime\prime}\in B^{**}$ and $(b^{\prime\prime}_\alpha)_\alpha\subseteq B^{**}$ such that $b^{\prime\prime}_\alpha\stackrel{w^*} {\rightarrow}b^{\prime\prime}$ in $B^{**}$. Then we have the following equality
$$ \lim_\alpha\langle  \pi_\ell^{***}(a^{\prime\prime},b^{\prime\prime}_\alpha),b^\prime\rangle
=\lim_\alpha\langle  \pi_\ell^{***}(a^{\prime\prime},b^{\prime\prime}_\alpha),x^\prime a \rangle $$$$=\lim_\alpha\langle  a\pi_\ell^{***}(a^{\prime\prime},b^{\prime\prime}_\alpha),x^\prime \rangle =
\langle  \pi_\ell^{***}(aa^{\prime\prime},b^{\prime\prime}),x^\prime \rangle =
\langle  \pi_\ell^{***}(a^{\prime\prime},b^{\prime\prime}),b^\prime \rangle .$$

It follows that $a^{\prime\prime}\in Z_{B^{**}}(A^{**}).$\end{proof}

\vspace{0.5cm}

In the preceding corollary, if we take $B=A$, then we have the following conclusion
$$\{a^{\prime\prime}\in A^{**}:~Aa^{\prime\prime}\subseteq A\}\subseteq Z_{1}(A^{**}).$$\\

\noindent{\it{\bf Theorem 4-3.}} Let $B$ be a Banach left $A-module$ and suppose that $b^\prime\in wap_\ell(B)$. Let  $a^{\prime\prime}\in A^{**}$ and
$(a_{\alpha})_{\alpha}\subseteq A$ such that  $a_{\alpha} \stackrel{w^*} {\rightarrow}a^{\prime\prime}$ in $A^{**}$. Then we have
$$\pi^*_\ell(b^\prime,a_\alpha) \stackrel{w} {\rightarrow}\pi_\ell^{****}(b^\prime,a^{\prime\prime}).$$
\begin{proof}
Assume that $b^{\prime\prime}\in B^{**}$. Then we have the following equality

$$\langle  \pi_\ell^{****}(b^\prime,a^{\prime\prime}),b^{\prime\prime}\rangle =
\langle  \pi_\ell^{***}(a^{\prime\prime},b^{\prime\prime}),b^\prime\rangle =
\lim_\alpha\langle  \pi_\ell^{***}(a_{\alpha},b^{\prime\prime}),b^\prime\rangle
$$$$=\lim_\alpha\langle  b^{\prime\prime},\pi_\ell^{*}(b^\prime, a_{\alpha})\rangle .$$
Now suppose that $(b_{\beta}^{\prime\prime})_{\beta}\subseteq B^{**}$ such that  $b^{\prime\prime}_{\beta} \stackrel{w^*} {\rightarrow}b^{\prime\prime}$. Since $b^\prime\in wap_\ell(B)$, we have
$$\langle  \pi_\ell^{****}(b^\prime,a^{\prime\prime}),b_\beta^{\prime\prime}\rangle =
\langle  \pi_\ell^{***}(a^{\prime\prime},b_\beta^{\prime\prime}),b^\prime\rangle \rightarrow
\langle  \pi_\ell^{***}(a^{\prime\prime},b^{\prime\prime}),b^\prime\rangle
$$$$=\langle  \pi_\ell^{****}(b^\prime,a^{\prime\prime}),b^{\prime\prime}\rangle .$$
Thus $\pi_\ell^{****}(b^\prime,a^{\prime\prime})\in (B^{**},weak^*)^*=B^*$.
So we conclude that
$$\pi^*_\ell(b^\prime,a_\alpha) \stackrel{w} {\rightarrow}\pi_\ell^{****}(b^\prime,a^{\prime\prime})~in ~ B^{**}.$$\\
\end{proof}

In the preceding corollary, if we take $B=A$, then we obtain the following result.\\
Suppose that $a^\prime\in wap(A)$ and $a^{\prime\prime}\in A^{**}$ such that
$a_{\alpha} \stackrel{w^*} {\rightarrow} a^{\prime\prime}$ where\\
$(a_{\alpha})_{\alpha}\subseteq A$. Then we have $a^\prime a_{\alpha} \stackrel{w}
{\rightarrow} a^\prime a^{\prime\prime}$.\\\\

\noindent{\it{\bf Theorem 4-4.}} Let $B$ be a Banach left $A-module$ and it has a $BLAI$ $(e_{\alpha})_{\alpha}\subseteq A$. Suppose that $b^\prime\in wap_\ell(B)$. Then we have
$$\pi^*_\ell(b^\prime,e_\alpha) \stackrel{w} {\rightarrow}b^\prime.$$
\begin{proof}
Let $b^{\prime\prime}\in B^{**}$ and $(b_\beta)_\beta\subseteq B$ such that $b_\beta\stackrel{w^*} {\rightarrow}b^{\prime\prime}$ in $B^{**}$. Then for every $b^\prime\in wap_\ell(B)$, we have the following equality
$$\lim_\alpha\langle  b^{\prime\prime}, \pi^*_\ell(b^\prime,e_\alpha)\rangle =
\lim_\alpha\langle  \pi^{****}_\ell( b^\prime,e_\alpha),b^{\prime\prime}\rangle $$$$=
\lim_\alpha\langle   b^\prime,\pi^{***}_\ell(e_\alpha,b^{\prime\prime})\rangle =
\lim_\alpha\langle   \pi^{***}_\ell(e_\alpha,b^{\prime\prime}),b^\prime\rangle $$$$=
\lim_\alpha\lim_\beta\langle   \pi_\ell(e_\alpha,b_\beta),b^\prime\rangle =
\lim_\beta\lim_\alpha\langle   \pi_\ell(e_\alpha,b_\beta),b^\prime\rangle $$$$=
\lim_\beta\langle   b_\beta,b^\prime\rangle =\langle   b^{\prime\prime},b^\prime\rangle .$$
It follows that
$$\pi^*_\ell(b^\prime,e_\alpha) \stackrel{w} {\rightarrow}b^\prime.$$
\end{proof}
\vspace{0.5cm}

\noindent{\it{\bf Corollary 4-5.}} Let $B$ be a Banach left $A-module$ and it has a $BLAI$ $(e_{\alpha})_{\alpha}\subseteq A$. Suppose that $ wap_\ell(B)=B^*$. Then $B^*$ factors on the left with respect to $A$.\\\\

\noindent{\it{\bf Corollary 4-6.}} Let $A$ be an Arens regular Banach algebra with $LBAI$. Then $A^*$ factors on the left.\\\\

\noindent{\it{\bf Example 4-7.}} i) Let $G$ be finite group. Then we have the following equality
$$M(G)^*L^1(G)=M(G)^*~and ~L^\infty (G)L^1(G)=L^\infty (G).$$
ii) Consider the Banach algebra $(\ell^1,.)$ that is Arens regular Banach algebra with unit element. Then we have $\ell^\infty .\ell^1=\ell^\infty$.\\
iii) Let $\mathcal{K}(E)$ be a compact operators from $E$ into $E$. Let  $E=\ell^p$ with $p\in (1,\infty)$. Then by using [6, Theorem 2.6.23] and [19 , Example 2.3.18] , $\mathcal{K}(\ell^p)$ is Arens regular and amenable, respectively. Since $\mathcal{K}(\ell^p)$ is amenable, by using [19 , Proposition 2.2.1], $\mathcal{K}(\ell^p)$ has a BAI. Therefore by preceding corollary, $\mathcal{K}(\ell^p)^*$ factors on the left.\\\\

\noindent{\it{\bf Theorem 4-8.}} Let $B$ be a Banach left $A-module$ and $B^*$ factors on the left with respect to $A$. If $AA^{**}\subseteq  Z_{B^{**}}(A^{**})$, then  $Z_{B^{**}}(A^{**})=A^{**}$.
\begin{proof}
Let $b^{\prime\prime}\in B^{**}$ and $(b^{\prime\prime}_\alpha)_\alpha\subseteq B^{**}$ such that $b^{\prime\prime}_\alpha\stackrel{w^*} {\rightarrow}b^{\prime\prime}$ in $B^{**}$. Suppose that $a^{\prime\prime}\in A^{**}$. Since $B^*$ factors on the left with respect to $A$,  for every $b^\prime\in B^*$, there are $x^\prime\in B^*$ and $a\in A$ such that $b^\prime=x^\prime a$. Since $aa^{\prime\prime}\in Z_{B^{**}}(A^{**})$, we have
$$\langle  \pi^{***}_\ell(a^{\prime\prime},b_\alpha^{\prime\prime}),b^\prime \rangle =
\langle  \pi^{***}_\ell(a^{\prime\prime},b_\alpha^{\prime\prime}),x^\prime a\rangle $$$$=
\langle  a\pi^{***}_\ell(a^{\prime\prime},b_\alpha^{\prime\prime}),x^\prime \rangle =
\langle  \pi^{***}_\ell(aa^{\prime\prime},b_\alpha^{\prime\prime}),x^\prime \rangle $$$$
\rightarrow \langle  \pi^{***}_\ell(aa^{\prime\prime},b^{\prime\prime}),x^\prime \rangle
=\langle  \pi^{***}_\ell(a^{\prime\prime},b^{\prime\prime}),b^\prime \rangle .$$
It follows that
$$\pi^{***}_\ell(a^{\prime\prime},b_\alpha^{\prime\prime})\stackrel{w^*} {\rightarrow}
\langle  \pi^{***}_\ell(a^{\prime\prime},b^{\prime\prime}),$$
and so $a^{\prime\prime}\in Z_{B^{**}}(A^{**})$.
\end{proof}
\vspace{0.5cm}

\noindent{\it{\bf Corollary 4-9.}} Let $A$ be a Banach algebra and $A^*$ factors on the left. If $AA^{**}\subseteq Z_1(A^{**})$, then $A$ is Arens regular.\\\\

\noindent{\it{\bf Theorem 4-10.}} Let $B$ be a Banach left $A-module$ and $B^{**}$ has a $LBAI$ with respect to $A^{**}$. Then $B^{**}$ has a left unit with respect to $A^{**}$.
\begin{proof}
Assume that $(e^{\prime\prime}_{\alpha})_{\alpha}\subseteq A^{**}$ is a $LBAI$ for $B^{**}$. By passing to a subnet, we may suppose that there is  $e^{\prime\prime}\in A^{**}$ such that $e^{\prime\prime}_{\alpha} \stackrel{w^*} {\rightarrow}e^{\prime\prime}$ in $A^{**}$. Then for every $b^{\prime\prime}\in B^{**}$ and $b^\prime\in B^*$, we have

$$\langle  \pi_\ell^{***}(e^{\prime\prime},b^{\prime\prime}),b^\prime\rangle =
\langle  e^{\prime\prime},\pi_\ell^{**}(b^{\prime\prime},b^\prime)\rangle =
\lim_\alpha \langle  e_\alpha^{\prime\prime},\pi_\ell^{**}(b^{\prime\prime},b^\prime)\rangle $$$$=
\lim_\alpha \langle  \pi_\ell^{***}(e_\alpha^{\prime\prime},b^{\prime\prime}),b^\prime\rangle =
\langle  b^{\prime\prime},b^\prime\rangle .$$
It follows that $\pi_\ell^{***}(e^{\prime\prime},b^{\prime\prime})=b^{\prime\prime}$.\end{proof}
\vspace{0.5cm}

\noindent{\it{\bf Corollary 4-11.}} Let $A$ be a Banach algebra and $A^{**}$ has a $LBAI$. Then $A^{**}$ has a left unit with respect to the first Arens product.\\\\

\noindent{\it{\bf Theorem 4-12.}} Let $B$ be a Banach left $A-module$ and it has a $LBAI$ with respect to $A$. Then we have the following assertions.
\begin{enumerate}
\item  $B^*$ factors on the left with respect to $A$ if and only if for each $b^\prime\in B^*$, we have $\pi^*_\ell(b^\prime,e_\alpha) \stackrel{w} {\rightarrow}b^\prime$ in $B^*$.
\item $B$ factors on the left with respect to $A$ if and only if for each $b\in B$, we have $\pi^*_\ell(b,e_\alpha) \stackrel{w} {\rightarrow}b$ in $B$.
\end{enumerate}
\begin{proof}\begin{enumerate}
\item Assume that  $B^*$ factors on the left with respect to $A$. Then  for every $b^\prime\in B^*$, there are $x^\prime\in B^*$ and $a\in A$ such that $b^\prime=x^\prime a$. Then for every $b^{\prime\prime}\in B^{**}$, we have
$$\langle  b^{\prime\prime}, \pi^*_\ell(b^\prime,e_\alpha)\rangle = \langle  b^{\prime\prime}, \pi^*_\ell(x^\prime a,e_\alpha)\rangle =\langle  b^{\prime\prime}, \pi^*_\ell(x^\prime ,ae_\alpha)\rangle $$$$=
\langle  \pi^{**}_\ell(b^{\prime\prime}, x^\prime), ae_\alpha)\rangle \rightarrow \langle  \pi^{**}_\ell(b^{\prime\prime}, x^\prime), a\rangle $$$$=
\langle  b^{\prime\prime},b^\prime\rangle .$$
It follows that
 $\pi^*_\ell(b^\prime,e_\alpha) \stackrel{w} {\rightarrow}b^\prime$ in $B^*$.
The converse by Cohen factorization theorem hold.

\item It is similar to the preceding proof.
\end{enumerate}\end{proof}
\vspace{0.5cm}

In the preceding theorem, if we take $B=A$, we obtain Lemma 2.1 from [14].\\\\

\noindent{\it{\bf Theorem 4-13.}} Let $B$ be a Banach left $A-module$ and $A$ has a $LBAI$ $(e_{\alpha})_{\alpha}\subseteq A$ such that $e_{\alpha} \stackrel{w^*} {\rightarrow}e^{\prime\prime}$ in $A^{**}$ where $e^{\prime\prime}$ is a left unit for $A^{**}$. Suppose that $Z^t_{e^{\prime\prime}}(B^{**})=B^{**}$. {Then, $B$ factors on the right with respect to $A$ if and only if
$e^{\prime\prime}$ is a left unit for $B^{**}$.
\begin{proof}
Assume that  $B$ factors on the right with respect to $A$. Then  for every $b\in B$, there are $x\in B$ and $a\in A$ such that $b= ax$. Then for every $b^\prime\in B^*$, we have
$$\langle  \pi^*_\ell(b^\prime,e_\alpha),b\rangle =\langle  b^\prime,\pi_\ell(e_\alpha,b)\rangle =\langle  \pi^{***}_\ell(e_\alpha,b),b^\prime\rangle $$$$=
\langle  \pi^{***}_\ell(e_\alpha,ax),b^\prime\rangle =
\langle  \pi^{***}_\ell(e_\alpha a,x),b^\prime\rangle $$$$=\langle  e_\alpha a,\pi^{**}_\ell(x,b^\prime)\rangle =
\langle  \pi^{**}_\ell(x,b^\prime),e_\alpha a\rangle $$$$\rightarrow \langle   \pi^{**}_\ell(x,b^\prime),a\rangle =\langle  b^\prime,b\rangle .$$
It follows that  $\pi^*_\ell(b^\prime,e_\alpha) \stackrel{w^*} {\rightarrow}b^\prime$ in $B^*$.
Let $b^{\prime\prime}\in B^{**}$ and $(b_\beta)_\beta\subseteq B$ such that $b_\beta\stackrel{w^*} {\rightarrow}b^{\prime\prime}$ in $B^{**}$. Since $Z^t_{e^{\prime\prime}}(B^{**})=B^{**}$, for every $b^\prime\in B^*$, we have the following equality
$$\langle  \pi_\ell^{***}(e^{\prime\prime},b^{\prime\prime}),b^\prime\rangle =
\lim_\alpha\lim_\beta\langle  b^\prime, \pi_\ell(e_\alpha,b_\beta)\rangle $$$$=
\lim_\beta\lim_\alpha\langle  b^\prime, \pi_\ell(e_\alpha,b_\beta)\rangle =
\lim_\beta\langle  b^\prime,b_\beta\rangle $$$$
=\langle  b^{\prime\prime},b^\prime\rangle .$$
It follows that $\pi_\ell^{***}(e^{\prime\prime},b^{\prime\prime})=b^{\prime\prime}$, and so $e^{\prime\prime}$ is a left unit for $B^{**}$.\\
Conversely, let $e^{\prime\prime}$ be a left unit for $B^{**}$ and suppose that $b\in B$. Thren for every $b^\prime\in B^*$, we have
$$\langle  b^\prime,\pi(e_\alpha,b)\rangle =\langle  \pi^{***}(e_\alpha,b),b^\prime\rangle =\langle  e_\alpha,\pi^{**}(b,b^\prime)\rangle
=\langle  \pi^{**}(b,b^\prime),e_\alpha\rangle $$$$=
\langle  e^{\prime\prime},\pi^{**}(b,b^\prime)\rangle =\langle  \pi^{***}(e^{\prime\prime},b),b^\prime\rangle =\langle  b^\prime,b\rangle .$$
Then we have $\pi^*_\ell(b^\prime,e_\alpha) \stackrel{w} {\rightarrow}b^\prime$ in $B^*$, and so by Cohen factorization theorem we are done.

\end{proof}
\vspace{0.5cm}

\noindent{\it{\bf Corollary 4-14.}} Let $B$ be a Banach left $A-module$ and $A$ has a $LBAI$ $(e_{\alpha})_{\alpha}\subseteq A$ such that $e_{\alpha} \stackrel{w^*} {\rightarrow}e^{\prime\prime}$ in $A^{**}$ where $e^{\prime\prime}$ is a left unit for $A^{**}$. Suppose that $Z^t_{e^{\prime\prime}}(B^{**})=B^{**}$. Then
$\pi^*_\ell(b^\prime,e_\alpha) \stackrel{w} {\rightarrow}b^\prime$ in $B^*$ if and only if $e^{\prime\prime}$ is a left unit for $B^{**}$.\\\\

\noindent{\it{\bf Problem.}} Let $B$ be a Banach  $A-bimodule$. Then\\
i) If $B$ factors, when $B^*$ factors?\\
ii) If $B$ is separable, dose $B$ necessarily factor on the one side?\\\\\\

\begin{center}
\textbf{{ 5. Involution $\ast-$algebra and Arens regularity of module actions}}
\end{center}
\vspace{0.5cm}

\noindent {\it{\bf Definition 5-1.}} Let $B$ be a Banach left $A-module$ and let $(a_\alpha)_\alpha\subseteq A$ has $weak^*$ limit in $A^{**}$. We say that $a_\alpha)_\alpha$ is  left regular with respect to $B$, if for every $(b_\beta)_\beta\subseteq B$,
$$w^*-\lim_\alpha w^*-\lim_\beta a_\alpha b_\beta =w^*-\lim_\beta w^*-\lim_\alpha a_\alpha b_\beta ,$$
where $(b_\beta)_\beta$ has $weak^*$ limit in $B^{**}$.\\
The definition of the right regular is similar. For a Banach $A-bimodule$ $B$, if $(a_\alpha)_\alpha\subseteq A$ is left and right regular with respect to $B$, we say that $(a_\alpha)_\alpha$ is  regular with respect to $B$. If $(a_\alpha)_\alpha$ is  left or right regular with respect to $A$, we write $(a_\alpha)_\alpha$ is  left or right regular, respectively.\\
\vspace{0.5cm}

\noindent {\it{\bf Example 5-2.}}\\ i) Let $B$ be a right Banach $A-module$ and $a^{\prime\prime}\in A^{**}$. Suppose that $Z^r_{a^{\prime\prime}}(B^{**})=B^{**}$ and $w^*-\lim_\alpha=a^{\prime\prime}$ where $(a_\alpha)_\alpha\subseteq A$. Then  $(a_\alpha)_\alpha$ is right regular with respect to $B$.\\
ii) Let $A$ be a Banach algebra and $(a_\alpha)_\alpha\subseteq A$ $weak^*$ convergence to some point of $Z_1(A^{**})$. Then it is clear that $(a_\alpha)_\alpha$ is regular.\\
\vspace{0.5cm}

\noindent {\it{\bf Theorem 5-3.}} i)  Let $B$ be a Banach left (resp. right) $A-module$ and suppose that $(e_{\alpha})_{\alpha}\subseteq A$ is a $LBAI$ for $B$. If $(e_{\alpha})_{\alpha}\subseteq A$ is left (resp. right) regular with respect to $B$, then $B^*$ factors on the left (resp. right).\\
ii) Let $B$ be a Banach left (resp. right) $A-module$ and suppose that  $B^{**}$ has  a left  (resp. right) unit as $A^{**}$-module. If $B^*A$ (resp. $AB^*$) is closed subspace of $B^*$, then  $B^*$ factors on the left (resp. right).
\begin{proof}  Let $e^{\prime\prime}\in A^{**}$ and  $e_{\alpha} \stackrel{w^*} {\rightarrow}e^{\prime\prime}$.  Assume that  $(b_\beta)_\beta\subseteq B$ such that $b_\beta \stackrel{w^*} {\rightarrow}b^{\prime\prime}$ in $B^{**}$. Then
$$e^{\prime\prime}b^{\prime\prime}=w^*-\lim_\alpha e_\alpha b^{\prime\prime}=w^*-\lim_\alpha w^* -\lim_\beta e_\alpha b_\beta=w^*-\lim_\beta w^* -\lim_\alpha e_\alpha b_\beta$$
$$=w^*-\lim_\beta b_\beta=b^{\prime\prime}.$$
Thus for every $b^\prime\in B^*$, we have
$$\langle  b^{\prime\prime},b^\prime\rangle  =\langle  e^{\prime\prime}b^{\prime\prime},b^\prime\rangle  =\lim_\alpha\langle  e_\alpha b^{\prime\prime},b^\prime\rangle  =\lim_\alpha\langle   b^{\prime\prime},b^\prime e_\alpha\rangle  .$$
It follows that  $b^\prime e_\alpha \stackrel{w^*} {\rightarrow}b^\prime$, and so by Cohen factorization theorem, we are done.\\
ii) Assume that $B^*A\neq B^*$. Let $e^{\prime\prime}\in A^{**}$ be a left unit element for $B^{**}$ and suppose that there is a net  $(e_{\alpha})_{\alpha}\subseteq A$ such that  $e_{\alpha} \stackrel{w^*} {\rightarrow}e^{\prime\prime}$.   By Hahn Banach theorem take $0\neq  b^{\prime\prime}\in B^{**}$ such that $ \langle  b^{\prime\prime}, B^*A\rangle  =0$. Then for every $b^\prime\in B^*$, we have
$$\langle  b^{\prime\prime},b^\prime\rangle  =\langle  e^{\prime\prime}b^{\prime\prime},b^\prime\rangle  =\lim_\alpha\langle  e_\alpha b^{\prime\prime},b^\prime\rangle  =\lim_\alpha\langle   b^{\prime\prime},b^\prime e_\alpha\rangle  .$$
That is contradiction. Thus $B^*A=B^*$.\\
Proof of the next part is similar.\end{proof}

\vspace{0.5cm}

\noindent {\it{\bf Definition 5-4.}}   Let $B$ be a Banach  $A-bimodule$ and suppose that $A$ is a Banach $\ast-involution$ algebra. We say that $B$ is a Banach $\ast-involution$ algebra as $A-module$, if the mapping $b\rightarrow b^*$ from $B$ into $B$ satisfies in the following conditions
$$(ab)^*=b^*a^*,~~(ba)^*=a^*b^*,~~(\lambda b)^*=\bar{\lambda}b^*,~~(b^*)^*=b~ ,~~ \parallel b^*\parallel,$$
for all $a\in A$, $b\in B$ and $\lambda\in \mathbb{C}$.\\

\vspace{0.5cm}

\noindent {\it{\bf Theorem 5-5.}}   Let $B$ be a Banach  $A-bimodule$ and suppose that $(e_{\alpha})_{\alpha}\subseteq A$ is a $LBAI$ for $B$. Let  $(e_{\alpha})_{\alpha}\subseteq A$ be a left  regular with respect to $B$. If  $B^{**}$ is a Banach $\ast-involution$ algebra as $A^{**}$-module, then  $B^{**}$ is unital  as $A^{**}$-module and $B^*$ factors on the both side.
\begin{proof} By using Theorem 5-3, $B^*$ factors on the left and without loss generally,  there is a left unit for $B^{**}$ as $e^{\prime\prime}\in A^{**}$  such that   $e_{\alpha} \stackrel{w^*} {\rightarrow}e^{\prime\prime}$. Let $b^{\prime\prime}\in B^{**}$. Then
$$b^{\prime\prime}(e^{\prime\prime})^*=(e^{\prime\prime}(b^{\prime\prime})^*)^*=((b^{\prime\prime})^*)^*=b^{\prime\prime}.$$
Since $e^{\prime\prime}=(e^{\prime\prime})^*$, $B^{**}$ is unital. Thus it is similar to Theorem 5-3, $B^*$ factors on the right.\end{proof}

\vspace{0.5cm}

\noindent {\it{\bf Corollary 5-6.}} Let $B$ be a Banach  $A-bimodule$ and suppose that $(e_{\alpha})_{\alpha}\subseteq A$ is a $BAI$ for $B$. Then if $B^{**}$ is a Banach $\ast-involution$ algebra as $A^{**}$-module,  then $B^{**}$ is unital  as $A^{**}$-module and $B^*$ factors on the both side.\\

\vspace{0.5cm}

\noindent {\it{\bf Example 5-7.}}  Let $G$ be a locally compact group. Then, on $L^1(G)$ there is a natural involution $*$ defined by $f^*(x)=\Delta (x^{-1})\overline{f(x^{-1})}$, where $\Delta$ is modular function and $x\in G$.
By preceding  theorem if $L^1(G)^{**}$ is a Banach $\ast-involution$ algebra, then $L^1(G)^{**}$ is unital with respect to first Arens product, and so $LUC(G)=L^\infty (G)$. It follows that $G$ is discrete.\\

\vspace{0.5cm}

\noindent {\it{\bf Problem .}}   Let $G$ be a locally compact group. We know that $M(G)$ by $\mu^*(f)=\overline{\int\overline{f(x^{-1})}d\mu}$ for all $f\in C_0(G)$ is $\ast$-involution algebra. When there is extension of it into $M(G)^{**}$ where $M(G)^{**}$ equipped first Arens product.\\\\

\bibliographystyle{amsplain}

\noindent Department of Mathematics, University of Mohghegh Ardabili, Ardabil, Iran\\
{\it Email address:} haghnejad@aut.ac.ir

\end{document}